\documentclass[onefignum,onetabnum]{siamart190516}

\usepackage[utf8]{inputenc}

\usepackage{amsmath}
\usepackage{amssymb}
\usepackage{bm}
\usepackage{soul}
\usepackage{multirow}
\usepackage{mathtools}

\usepackage{tikz}
\usetikzlibrary{arrows.meta}
\usetikzlibrary{calc}
\usetikzlibrary{decorations.markings, decorations.pathreplacing, angles,quotes}

\newtheorem{remark}{Remark}
\newtheorem{assumption}{Assumption}
\newtheorem{example}{Example}

\newcommand{\average}[1]{\overline{#1}}
\newcommand{\distributionsadaptive}{f}
\newcommand{\momentumadaptive}{m}
\newcommand{\detailsadaptive}{d}
\newcommand{\projectionoperator}{\mathbf{P}_{\vee}}
\newcommand{\predictionoperator}{\mathbf{P}_{\wedge}}
\newcommand{\reals}{\mathbb{R}}
\newcommand{\levelletter}{j}
\newcommand{\indexletter}{k}
\newcommand{\populationletter}{h}
\newcommand{\adaptiveroundbrackets}[1]{\left ( #1 \right )}
\newcommand{\adaptivecurlybrackets}[1]{\left \{ #1 \right \}}
\newcommand{\ratio}[2]{\dfrac{#1}{#2}}
\newcommand{\ratioobl}[2]{#1 / #2}
\newcommand{\velocitynumber}{q}
\newcommand{\minlevel}{\underline{J}}
\newcommand{\maxlevel}{\overline{J}}
\newcommand{\levelsnumber}{L}
\usepackage{scalerel,stackengine}
\stackMath
\newcommand\reallywidehat[1]{%
\savestack{\tmpbox}{\stretchto{%
  \scaleto{%
    \scalerel*[\widthof{\ensuremath{#1}}]{\kern-.6pt\bigwedge\kern-.6pt}%
    {\rule[-\textheight/2]{1ex}{\textheight}}
  }{\textheight}%
}{0.5ex}}%
\stackon[2pt]{#1}{\tmpbox}%
}
\newcommand\reallywidedoublehat[1]{%
\savestack{\tmpbox}{\stretchto{%
  \scaleto{%
    \scalerel*[\widthof{\ensuremath{#1}}]{\kern-.6pt\bigwedge\kern-.6pt}%
    {\rule[-\textheight/2]{1ex}{\textheight}}
  }{\textheight}%
}{0.5ex}}%
\stackon[-0.3mm]{\stackon[1pt]{#1}{\tmpbox}}{\tmpbox}%
}

\newcommand{\superscript}[2]{#1^{#2}}
\newcommand{\subscript}[2]{#1_{#2}}
\newcommand{\predicted}[3]{\reallywidehat{#1} \subscript{\superscript{\vphantom{#1}}{#2}}{#3}}
\newcommand{\reconstructed}[3]{\reallywidedoublehat{#1} \subscript{\superscript{\vphantom{#1}}{#2}}{#3}}
\newcommand{\definitionequality}{:=}
\newcommand{\setofindices}{\nabla}

\newcommand{\vectorial}[1]{\bm{#1}}
\newcommand{\multiresolutiontransform}{\mathcal{M_R}}
\newcommand{\physicaltree}[1]{R(#1)}
\newcommand{\leaves}[1]{L(#1)}
\newcommand{\physicalleaves}[1]{S(#1)}
\newcommand{\predictionstencilwidth}{\gamma}
\newcommand{\gradingoperator}{\mathcal{G}}
\newcommand{\thresholdletter}{\epsilon}
\newcommand{\thresholdoperator}[1]{\mathcal{T}_{\vectorial{\thresholdletter}}(#1)}
\newcommand{\norm}[1]{\left \lVert #1 \right \rVert}
\newcommand{\enlargingoperator}{\mathcal{H}_{\bm{\thresholdletter}}}
\newcommand{\adaptiveschemeoperator}{\bm{L}_{A}}
\newcommand{\referenceschemeoperator}{\bm{L}}
\newcommand{\finaltime}{T}
\newcommand{\latticevelocity}{\lambda}
\newcommand{\velocityletter}{v}
\newcommand{\normalizedvelocityletter}{w}
\newcommand{\operatorial}[1]{\bm{#1}}
\newcommand{\predictionstencildepth}{\gamma}
\newcommand{\scheme}[2]{$\text{D}_{#1}\text{Q}_{#2}$}
\newcommand{\schemevec}[3]{$\text{D}_{#1}\text{Q}_{#2}^{#3}$}

\newcommand{\myrightleftarrows}[1]{\mathrel{\substack{\xrightarrow{#1} \\[-.9ex] \xleftarrow{#1}}}}

\newcommand{\additionRefOne}[1]{{#1}}
\newcommand{\additionRefTwo}[1]{{#1}}
\newcommand{\additionRefTwoSecondRound}[1]{#1}
\newcommand{\additionFree}[1]{{#1}}

\usepackage{soul}
\usepackage{calc}
\newsavebox\CBox
\newcommand\hcancel[2][0.5pt]{%
  \ifmmode\sbox\CBox{$#2$}\else\sbox\CBox{#2}\fi%
  \makebox[0pt][l]{\usebox\CBox}%
  \rule[0.5\ht\CBox-#1/2]{\wd\CBox}{#1}}
\newcommand{\strikeformula}[1]{\textcolor{orange}{\hcancel{#1}}}

\title{Multiresolution-based mesh adaptation and error control for lattice Boltzmann methods with applications to hyperbolic conservation laws}

\author{Thomas Bellotti\thanks{CMAP, CNRS, Ecole polytechnique, Institut Polytechnique de Paris, 91128 Palaiseau Cedex, France. \email{thomas.bellotti@polytechnique.edu}} \and Lo\"{I}c Gouarin\footnotemark[1] \and Benjamin Graille\thanks{Institut de Mathématiques d'Orsay, Université Paris-Saclay, 91405 Orsay Cedex, France.} \and Marc Massot\footnotemark[1]}

\begin{document}

\maketitle

\begin{abstract}
    Lattice Boltzmann Methods (LBM) stand out for their simplicity and computational efficiency while offering the possibility of simulating complex phenomena. 
    While they are optimal for Cartesian meshes, adapted meshes have traditionally been a stumbling block since it is difficult to predict the right physics through various levels of meshes.
     In this work, we design a class of fully adaptive LBM methods with dynamic mesh adaptation and error control relying on multiresolution analysis. This wavelet-based approach allows to 
     adapt the mesh based on the regularity of the solution and leads to a very efficient compression of the solution without loosing its quality and 
     with the preservation of the properties of the original LBM method on the finest grid. This yields a general approach for a large spectrum of schemes and allows precise error bounds, 
     without the need for deep modifications on the reference scheme. An error analysis is proposed. 
     For the purpose of assessing the approach, we conduct a series of test-cases for various schemes and scalar and systems of conservation laws, 
     where solutions with shocks are to be found and local mesh adaptation is especially relevant. Theoretical estimates are retrieved while a reduced memory footprint is observed. 
     It paves the way to an implementation in a multi-dimensional framework and high computational efficiency of the method for both parabolic and hyperbolic equations, which is the subject of  a companion paper.
\end{abstract}



%
%
%
%
%
%

\section{Introduction}



A wide class of systems representing  various complex phenomena across different disciplines (fluid mechanics, 
combustion, atmospheric sciences, plasma physics or biomedical engineering, see \cite{duarte2013,duarte2014,duarte2015} and references therein for a few examples)  are modeled through PDEs, the solution of which
can involve dynamically moving fronts, usually very localized in space. Among these PDEs, one can find the fluid dynamics Euler equations, and more generally 
hyperbolic systems of conservation laws, where shock wave solutions are to be found. 
For such solutions,  we need a good level of spatial detail where steep variations occur, whereas one can accept a coarse space discretization where large \emph{plateaux} are present.
An effective way of reducing the overall cost of a numerical solvers consists in devising a strategy to dynamically adapt the spatial discretization to the solution as time advances, 
aiming at performing less operations and limiting the memory footprint, while preserving a proper resolution.  
Once a discretization is chosen, there exists several strategies for mesh adaptation in terms of both data-structure and refinement strategy. 
Such strategies can make a crucial difference in terms of time-to-solution and allow scientists to strongly reduce computational cost or reach the solution of large 3D problems on standard machines.
In terms of data-structure, one can either choose patch-based/block-based \cite{ray2007using, narechania2017block} refinement of cell-based/multiresolution \cite{burstedde2011p4est, cohen2003}, the first one being easier to parallelize, while the second is more optimal in terms of compression rate. Beyond such a choice, adapting the mesh relies on a choice of refinement criteria. There are several possibilities ranging from feature-based \cite{guittet2015stable}, discretization errors \cite{naddei2019comparison} as well as a posteriori estimates \cite{wu1990error, alauzet2003transient}, Richardson error evaluation \cite{berger1984adaptive2}, goal-oriented criteria relying on adjoint evaluation \cite{narechania2017block} or even optimal sparse sensing \cite{foti2020adaptive} to cite a few. Our purpose here is to tackle unsteady problems of hyperbolic and parabolic type, where the mesh is going to evolve dynamically and we aim at providing error control on the solution, so that we focus on the multiresolution approach. The wavelet-based approach coupled with the Harten heuristics leads in practice to solution error control based on the regularity of the solution. Even if optimality for steady problem is hard for such approches, it is well adapted to the typical problems solved by LBM methods.

The discretization of the original PDEs can be conducted relying on several methods: 
here we focus on Lattice Boltzmann schemes (LBM - Lattice Boltzmann Methods), a class of wide-spread numerical methods to approximate models which can have spatially non-homogeneous solutions. Despite being present in the community since the end of the eighties, they have gained a lot of attention in the last decade due to the evolution of the computer architectures and body of literature on their mathematical analysis. Mesh-adaptation for LBM has been a stumbling block for quite some time even if interesting pieces of solution have been provided.
The key issue is related to the difficulty of predicting the right physics through various levels of meshes. This can also lead to delicate transmission conditions for acoustic waves and has been a relatively hot topic in the field.

The purpose of the present contribution is to design of a new numerical strategy for LBM with dynamic mesh adaptation
 and error control based on multiresolution analysis.
The key issue is the ability to rely on the original LBM scheme on the finest mesh without alteration while still reaching a controlled level of accuracy on the compressed representation.
This contribution focuses on setting the fundamentals of the method, we concentrate on the one-dimensional framework and provide an error analysis. We conduct a numerical assessment on various test-cases for hyperbolic conservation laws (scalar and systems). We have chosen such a framework since it is a very representative example of localized fronts with default of regularity, where the MR can play a key role and where we can test a large variety of LBM schemes. The method yields a very reduced memory footprint while preserving a given level of accuracy. 
The proposed numerical strategy is versatile and can be extended in a straightforward manner to parabolic and hyperbolic systems in multi-dimensions, which is out of the scope of the present paper but is the subject of a companion contribution \cite{bellotti2021}.  Before entering the body of the contribution, let us describe the state of the art. 

%

\subsection{Lattice Boltzmann methods and mesh adaptation}

The lattice Boltzmann methods are relatively recent computational techniques for the numerical solution of PDEs, 
introduced at the end of the eighties by McNamara and Zanetti \cite{mcnamara1988} and by Higuera and Jimenez \cite{higuera1989}, and stemming from the 
 ``Lattice gas automata''.
The derivation of the method starts from Boltzmann equation, 
with a simplified collision kernel\footnote{Through the BGK approximation with single or multiple relaxation times.}, and relies on the  selection of a small set of discrete velocities compatible with a given fixed-step lattice.
This strategy is widely employed in many areas of computational mathematics, with special mention to the Computational Fluid Dynamics.
In this context, the method has been used to simulate the Navier-Stokes system at low Mach numbers \cite{lallemand2000} with more recent extensions to handle multi-phase problems (\cite{huang2015} for a review), along with systems of hyperbolic conservation laws 
\cite{graille2014}.
The advantages of the method are its dramatic simplicity\footnote{Since overall the strategy decomposes into a local collision and a stream along the characteristics of the discrete velocities} and the ease of parallelization.
Still, stability, consistency and convergence remain open topics.

To the best of our knowledge, LBM strategies on adapted grids have been only developed either on fixed grids, in the spirit of Filippova and H{\"a}nel \cite{filippova1998} and of many subsequent works, where more refined patches are placed according to an \emph{a priori} knowledge of the flow. Such fixed refinement zones also yield difficulties in aeroacoustics resolution 
related to the artificial transmission impedance of the refinement interface
\cite{gendre2017,feng2020, horstmann2018hybrid}.
Another strategy is to use an AMR approach \cite{berger1989} with some heuristics to determine the need for refinement in certain areas.
In this class, we find the work of Fakhari and Lee \cite{fakhari2014} using the magnitude of the vorticity and its derivatives as regularity indicator, while Eitel-Amor \emph{et al.} \cite{eitel2013} have employed a weighted vorticity and the energy difference with respect to a free flow solution.
Crouse \emph{et al.} \cite{crouse2003} have used the weighted magnitude of the divergence of the velocity field.
Finally, Wu and Shu \cite{wu2011} have considered the difference between solutions at successive time steps. 
Although these approaches have been certainly able to reduce the computational cost of the simulations, they still face several drawbacks which we summarize as follows:
\begin{itemize}
    \item Few available methods are time-adaptive: most of the time, one must construct a fixed non-uniform mesh according to some \emph{a priori} knowledge of the solution.
        The refinement interface can then induce spurious effects on the numerical simulations when fine-scale physics interfere with a coarse level of the mesh. An example of such a situation is the resolution of acoustic waves leading to purely numerical transmission defaults. 
    Consequently, this approach is intrinsically problem-dependent and non-optimal for unsteady solutions. 
    \item Most of the time the reference scheme has to be deeply manipulated at various levels of grid to preserve the macroscopic parameters of the system.
    \item One must devise a heuristic, the dynamic mesh adaptation will rely on. As a consequence, there is no control on the perturbation error by the mesh adaptation.
\end{itemize}

In this work, we propose a strategy to fill these gaps by introducing a time adaptive numerical approach, designed to work for any LBM scheme without need for manipulations, which guarantees a precise bound on the perturbation error.

\subsection{Multiresolution analysis}

Multiresolution analysis has proved to be a general tool to analyze the local regularity of a signal in a rigorous setting, based on its decomposition on a wavelet basis.
It has been introduced in the seminal works by Daubechies \cite{daubechies1988}, Mallat \cite{mallat1989} and Cohen \emph{et al.} \cite{cohen1992}.
The possibility of applying this mechanism to reduce the computational cost of a numerical method was studied a few years later by Harten \cite{harten1994, harten1995,bihari1997} in the context of Finite Volume methods for conservation laws.
The principle was to use multiresolution to reduce the number of computations to evaluate fluxes at the interfaces, claiming that they constitute the majority of the computational cost.
However, this approach still computes the solution on the full uniform mesh.
The possibilities offered by multiresolution had been further exploited by Cohen \emph{et al.} \cite{cohen2003} who, in the footsteps of Harten, have developed fully adaptive schemes with solutions updated only on the reduced grid.
Thus, multiresolution is not only a way of computing a large number of fluxes more cheaply, but also a manner to compute fewer of them.
Both these strategies ensure better time-performances than traditional approaches on uniform grids in addition to a precise control on the perturbation error, unlike most of the AMR techniques.
This strategy has been lately used to tackle various kinds of problems with Finite Volume methods.
We mention parabolic conservation laws by Roussel \emph{et al.} \cite{roussel2003}, the compressible Navier-Stokes equations in Bramkamp \emph{et al.} \cite{bramkamp2004}, the shallow water equations by Lamby \emph{et al.} \cite{lamby2005}, multi-component flows by Coquel \emph{et al.} \cite{coquel2006}, degenerate parabolic equations by Burger \emph{et al.} \cite{burger2008} and finally the Euler system with a local time-stepping technique again by Coquel \emph{et al.} \cite{coquel2010}.
Furthermore, this technique has been included in later works to address more complex problems, such as flames \cite{roussel2005,duartecandel13,duarte2014} or by coupling it with other numerical strategies: we mention the works of Duarte \emph{et al.} \cite{duarte2012,duarte2013,duarte2015} and N'Guessan \emph{et al.} \cite{nguessan2019}. 
We are not aware of the use of such procedure to conduct mesh adaptation and error control on LBM schemes.
We decided to adopt a volumetric vision for MR because since it is naturally conservative, even if a whole body of literature exists about point-wise multiresolution \cite{harten1993discrete,chiavassa2001,forster2016}.

\subsection{Paper organization}

We present the formalism of LBM schemes and Multiresolution analysis in section \ref{sec:LBMSchemes} and \ref{sec:AdaptiveMultiresolution};  the fundamentals of the proposed numerical strategy and the theoretical error analysis is presented in section \ref{sec:AdaptiveLBMMR}. Section \ref{sec:Verifications} is dedicated to numerical verifications on scalar and hyperbolic systems of conservation laws with a variety of LBM schemes, before concluding in section \ref{sec:Conclusions}.

\section{Lattice Boltzmann schemes}\label{sec:LBMSchemes}

Consider a uni-dimensional bounded domain $\Omega = [a, b]$ with $a < b$ and the maximum level of allowed refinement $\maxlevel \in \mathbb{N}$.
The domain is discretized by a partition of $2^{\maxlevel}$ \additionFree{cells} with measure $\Delta x_{\maxlevel} = 2^{-\maxlevel} (b-a)$ forming a collection $\mathcal{L}_{\maxlevel} \definitionequality (I_{\maxlevel, \indexletter})_{\indexletter  = 0, \dots, 2^{\maxlevel} - 1}$ called lattice.
Once we consider a function $f(t, x)$ of time and space, we define -- in a Finite Volume fashion -- its spatial averages on each cell
\begin{equation*}
    \average{F}_{\indexletter}(t) \approx \ratio{1}{|I_{\maxlevel, \indexletter}|} \int_{I_{\maxlevel, \indexletter}} f(t, x)\text{d}x, \qquad t \geq 0, \quad \indexletter = 0, \dots, 2^{\maxlevel} - 1.
\end{equation*}
Henceforth, every discretized quantity \additionRefTwo{with a bar} shall be interpreted as a mean value of an underlying integrable function over the cell it refers to.
In all the work, we consider a finite temporal horizon $t \in [0, \finaltime]$ with $\finaltime > 0$. The time is discretized, as we shall see in a moment, in equally spaced time steps with step $\Delta t > 0$.
Without loss of generality, we assume that $\finaltime$ has been chosen so that $N \definitionequality T/\Delta t \in \mathbb{N}$. Thus we indicate $t^n \definitionequality n\Delta t$ for $n = 0, \dots, N$ and $\average{F}^n \approx \average{F}_k(t^n)$ for $k = 0, \dots, 2^{\maxlevel} - 1$ and $n = 0, \dots, N$.
From now on $\lVert \cdot \rVert_{\ell^p}$ shall denote the weighted $\ell^p$ norm for $p \in [1, \infty]$ over $\reals^{2^{\maxlevel}}$ with weight $\Delta x_{\maxlevel}$.

\subsection{The d'Humières formalism}

We consider lattice Boltzmann schemes under the so-called d'Humières formalism \cite{dhumieres1992}.
Let $\latticevelocity > 0$ be a lattice velocity, so that the time step $\Delta t$ can be defined using the acoustic scaling\footnote{Still, the strategy of this work equally works for a parabolic scaling such as $\Delta t \sim (\Delta x_{\maxlevel})^2$ \cite{bellotti2021}.} $\Delta t = \Delta x_{\maxlevel} / \lambda$.
Moreover, consider a set of discrete velocities $\{\velocityletter^{\populationletter} \}_{\populationletter = 0, \dots, \velocitynumber - 1}$, where $\velocitynumber \in \mathbb{N}$, compatible with the lattice velocity $\latticevelocity$ in the sense that $\normalizedvelocityletter^{\populationletter} \definitionequality \velocityletter^{\populationletter} / \latticevelocity \in \mathbb{Z}$ for $\populationletter = 0, \dots, \velocitynumber - 1$.
At time $t^n$, we indicate with $\average{F}^{\populationletter, n}_{\indexletter}$  the average of the density of the population moving with velocity $\velocityletter^{\populationletter}$ on the cell $\indexletter = 0, \dots, 2^{\maxlevel} - 1$ for every $\populationletter = 0, \dots, \velocitynumber - 1$.
In order to recover the so-called moments, we consider an invertible matrix $\operatorial{M} \in \text{GL}_{\velocitynumber}(\reals)$ defining the change of variables to pass from the space of the population densities towards the space of moments by $(\average{M}^0, \dots,\average{M}^{\velocitynumber - 1})^T = \operatorial{M} (\average{ F}^0, \dots, \average{F}^{\velocitynumber - 1})^T$, where we do not indicate the space and the time coordinates because the change of basis is completely local.
\additionRefTwo{Observe that the change of basis $\operatorial{M}$ utilized here for the averages is exactly the same that one could use with lattice Boltzmann schemes where the discretized quantities are interpreted as point values. This comes from the linearity of the integral.}

The lattice Boltzmann scheme can be divided into two phases: collision phase and stream phase.
We indicate its action as
\begin{equation*}
    (\average{F}^{\populationletter, n+1}_{\indexletter})_{\substack{k = 0, \dots, 2^{\maxlevel} - 1 \\ \populationletter = 0, \dots, \velocitynumber - 1}} = \referenceschemeoperator (\average{F}^{\populationletter, n}_{\indexletter})_{\substack{k = 0, \dots, 2^{\maxlevel} - 1 \\ \populationletter = 0, \dots, \velocitynumber - 1}},
\end{equation*} 
for any $n = 0, \dots, N-1$.
\additionRefTwo{The link between this average formulation of the lattice Boltzmann schemes and the discrete Boltzmann equation is provided in the Supplementary material.}
The operator $\referenceschemeoperator$ is constructed as follows. 

\subsubsection{Collision phase}
Let us consider a cell $\indexletter = 0, \dots, 2^{\maxlevel} - 1$, then the collision phase is a local linear relaxation of the non-conserved moments towards their equilibrium, namely
\begin{align*}
    \average{M}_{\indexletter}^{\populationletter, n\star} = \average{M}_{\indexletter}^{\populationletter, n}, \qquad &\populationletter = 0, \dots, \velocitynumber_{\text{cons}} - 1, \\
    \average{M}_{\indexletter}^{\populationletter, n\star} = (1-s^{\populationletter}) \average{M}_{\indexletter} ^{\populationletter, n} + s^{\populationletter} M^{\populationletter, \text{eq}} (\average{M}_{\indexletter} ^{0, n}, \dots, \average{M}_{\indexletter} ^{\velocitynumber_{\text{cons}}-1, n} ), \qquad &\populationletter = \velocitynumber_{\text{cons}}, \dots, \velocitynumber - 1,
\end{align*}
where $\velocitynumber_{\text{cons}} < \velocitynumber$ is the number of conserved moments 
and where $s^{\populationletter}$ and $M^{\populationletter, \text{eq}}$ are respectively the relaxation parameter and the equilibrium of the $\populationletter^{\text{th}}$ moment, which is a non-linear function of the conserved moments. 
These quantities are set relying either on a Chapman-Enskog expansion or on the theory of equivalent equations introduced by Dubois \cite{dubois2009} in order to be consistent with the equations we want to solve.
\additionRefTwo{The relaxation parameters satisfy $s^{\populationletter} \in (0, 2]$, see \cite{dubois2009}. In the context of hyperbolic conservation laws, they influence the stability (\emph{i.e.} the numerical diffusion) of the method and also higher order terms in the resulting equivalent equations. The task of devising lattice Boltzmann schemes with good features tuning the relaxation parameters is a vast subject beyond the scope of this work, which highly depends on the considered method.}
\additionRefTwo{It is important to stress that the lattice Boltzmann method is used to solve a system of $\velocitynumber_{\text{cons}}$ conservation laws (PDEs), corresponding to the variables of interest, namely the conserved moments, by enlarging the size of the problem to $\velocitynumber > \velocitynumber_{\text{cons}}$. There is therefore some latitude on the choice of the initial condition (we take it at equilibrium), because the initial datum is known only on the conserved variables.}

\subsubsection{Stream phase}

Again on a cell $\indexletter = 0, \dots, 2^{\maxlevel} - 1$, the stream phase is given by
\begin{equation*}
    \average{F}^{\populationletter, n+1}_{\indexletter} = \average{F}^{\populationletter, n\star}_{\indexletter - \normalizedvelocityletter^{\populationletter}}, \qquad \populationletter = 0, \dots, \velocitynumber - 1.
\end{equation*}

As we shall be interested in considering all the populations together, in the sequel, the weighted $\ell^p$ norm are extended from $\reals^{2^{\maxlevel}}$ to $\reals^{\velocitynumber 2^{\maxlevel}}$ in the usual way.
\additionRefTwo{Moreover in what follows, uppercase letters, such as $\average{F}$ and $\average{M}$, indicate the solution of the reference scheme, namely the lattice Boltzmann scheme on the uniform lattice at maximum level $\maxlevel$.}
\additionRefTwo{The choice of employing a less known formulation based on averages as done by \cite{rohde2006} allows to easily enforce the conservation constraints when performing the mesh adaptation.}

\section{Adaptive multiresolution}\label{sec:AdaptiveMultiresolution}

Following the approach by \cite{cohen2003, duarte2011}, the starting point of the multiresolution analysis is to consider a hierarchy of $\levelsnumber + 1$ with $\levelsnumber \in \mathbb{N}$ nested uni-variate lattices $\mathcal{L}_{\levelletter}$ with $\levelletter = \maxlevel - \levelsnumber, \dots, \maxlevel$, given by
\begin{equation*}
    \mathcal{L}_{\levelletter} \definitionequality (I_{\levelletter, \indexletter})_{\indexletter = 0,\dots, N_{\levelletter} - 1}, \qquad \text{with} \quad I_{\levelletter, \indexletter} \definitionequality [(b-a)2^{-\levelletter}\indexletter + a, (b-a)2^{-\levelletter}(k+1)+a],
\end{equation*}
for $\levelletter = \minlevel, \dots, \maxlevel$, where we have set $\minlevel \definitionequality \maxlevel - \levelsnumber$ and $N_{\levelletter} \definitionequality 2^{\levelletter}$ for $\levelletter = \minlevel, \dots, \maxlevel$.
They form a sequence of progressively finer nested lattices as we observe that each cell $I_{\levelletter, \indexletter}$ (called ``parent'') includes its two ``children'' $I_{\levelletter + 1, 2\indexletter}$ and $I_{\levelletter + 1, 2\indexletter + 1}$ (called ``siblings'') rendering a tree-like structure.
As done before, given a function $f(t, x)$, we have to understand things in the following way
\begin{equation*}
    \average{\distributionsadaptive}_{\levelletter, \indexletter}^{n} \approx \frac{1}{|I_{\levelletter, \indexletter}|} \int_{I_{\levelletter, \indexletter}} f(t^n, x) \text{d}x,
\end{equation*}
for $n = 0, \dots, N$, $\levelletter = \minlevel, \dots, \maxlevel$ and $\indexletter = 0, \dots, N_{\levelletter} - 1$.
\additionRefTwo{We now use lowercase letters such as $\average{\distributionsadaptive}$ and $\average{m}$ to signify that we consider the solution of the adaptive lattice Boltzmann scheme which we are going to devise.}
In the remaining part of this Section, since time is of no importance, we do not mention the dependence of any quantity on it for the sake of clarity.

\subsection{Projection and prediction operator}

\begin{figure}[h]
        \begin{center}
            \begin{tikzpicture}[x=0.092cm, y=0.092cm]
                \newcount\levdis
                \levdis = 10
                
                \draw[Bar-Bar, very thick] (0, 0) -- (32, 0);
                \draw[-Bar, very thick]    (32, 0) -- (64, 0);
                \draw[-Bar, very thick]    (64, 0) -- (96, 0);
                \draw[-Bar, very thick]    (96, 0) -- (128, 0) node [right]  {$j $};

                \draw[Bar-Bar, very thick] (0, \levdis) -- (16, \levdis);
                \draw[-Bar, very thick]   (16, \levdis) -- (32, \levdis);
                \draw[-Bar, very thick]   (32, \levdis) -- (48, \levdis);
                \draw[-Bar, very thick]   (48, \levdis) -- (64, \levdis);
                \draw[-Bar, very thick]   (64, \levdis) -- (80, \levdis);
                \draw[-Bar, very thick]   (80, \levdis) -- (96, \levdis);
                \draw[-Bar, very thick]   (96, \levdis) -- (112, \levdis);
                \draw[-Bar, very thick]   (112, \levdis) -- (128, \levdis) node [right]  {$j + 1$};
                
                \draw[-stealth]  (8, 0.9*\levdis) -- (16-0.1*\levdis, 0.1*\levdis) node [midway, left]  {$\ratioobl{1}{2}$};
                \draw[-stealth] (24, 0.9*\levdis) -- (16+0.1*\levdis, 0.1*\levdis) node [midway, right] {$\ratioobl{1}{2}$};

                \draw[-stealth]  (40, 0.9*\levdis) -- (48-0.1*\levdis, 0.1*\levdis) node [midway, left]  {$\ratioobl{1}{2}$};
                \draw[-stealth]  (56, 0.9*\levdis) -- (48+0.1*\levdis, 0.1*\levdis) node [midway, right] {$\ratioobl{1}{2}$};
                
                \draw[-stealth]  (72, 0.9*\levdis) -- (80-0.1*\levdis, 0.1*\levdis) node [midway, left]  {$\ratioobl{1}{2}$};
                \draw[-stealth]  (88, 0.9*\levdis) -- (80+0.1*\levdis, 0.1*\levdis) node [midway, right] {$\ratioobl{1}{2}$};
                
                \draw[-stealth]  (104, 0.9*\levdis) -- (112-0.1*\levdis, 0.1*\levdis) node [midway, left]  {$\ratioobl{1}{2}$};
                \draw[-stealth]  (120, 0.9*\levdis) -- (112+0.1*\levdis, 0.1*\levdis) node [midway, right] {$\ratioobl{1}{2}$};
                
            \end{tikzpicture}
        \end{center}
        \caption[]{\label{fig:projection} Illustration of the action of the projection operator. The cell average on the cell at level $\levelletter$ is reconstructed by taking the average of the values on its two children at level $\levelletter + 1$.}
    \end{figure}
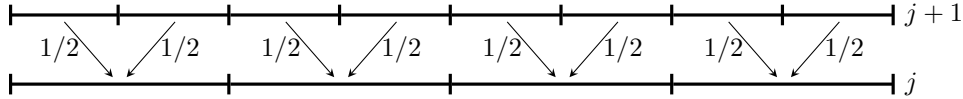

The projection and the prediction operators  allow one to navigate in the ladder of nested lattices up and down.
We start from the projection operator, which takes information at a certain level of resolution $j+1$ and transforms it into information on a coarser level $j$ as illustrated in Figure \ref{fig:projection}.

\begin{definition}[Projection operator]
    The projection operator $\projectionoperator : \reals^2 \to \reals$ is defined by
    \begin{equation*}
        \average{\distributionsadaptive}_{\levelletter, \indexletter}^{\populationletter} = \projectionoperator \adaptiveroundbrackets{\adaptiveroundbrackets{\average{\distributionsadaptive}^{\populationletter}_{\levelletter+1, 2\indexletter + \delta}}_{\delta = 0, 1}} = \ratio{1}{2} \adaptiveroundbrackets{\average{\distributionsadaptive}^{\populationletter}_{\levelletter+1, 2\indexletter } + \average{\distributionsadaptive}^{\populationletter}_{\levelletter+1, 2\indexletter + 1}},
    \end{equation*}
    for every $\populationletter = 0, \dots, \velocitynumber-1$, $\levelletter = \minlevel, \dots, \maxlevel$ and $\indexletter = 0, \dots, N_{j} - 1$.
\end{definition}

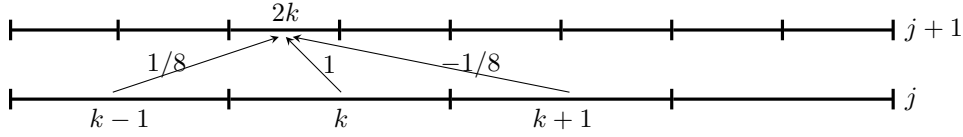
\begin{figure}[h]
        \begin{center}
            \begin{tikzpicture}[x=0.092cm, y=0.092cm]
                \newcount\levdis
                \levdis = 10
                
                \draw[Bar-Bar, very thick] (0, 0) -- (32, 0);
                \draw[-Bar, very thick]    (32, 0) -- (64, 0);
                \draw[-Bar, very thick]    (64, 0) -- (96, 0);
                \draw[-Bar, very thick]    (96, 0) -- (128, 0) node [right]  {$j $};

                \draw[Bar-Bar, very thick] (0, \levdis) -- (16, \levdis);
                \draw[-Bar, very thick]   (16, \levdis) -- (32, \levdis);
                \draw[-Bar, very thick]   (32, \levdis) -- (48, \levdis);
                \draw[-Bar, very thick]   (48, \levdis) -- (64, \levdis);
                \draw[-Bar, very thick]   (64, \levdis) -- (80, \levdis);
                \draw[-Bar, very thick]   (80, \levdis) -- (96, \levdis);
                \draw[-Bar, very thick]   (96, \levdis) -- (112, \levdis);
                \draw[-Bar, very thick]   (112, \levdis) -- (128, \levdis) node [right]  {$j + 1$};
                
                \draw[-stealth] (16-0.1*\levdis, 0.1*\levdis) -- (40-0.1*\levdis, 0.9*\levdis) node [midway, left]  {$\ratioobl{1}{8}$};
                \draw[-stealth] (48, 0.1*\levdis) -- (40, 0.9*\levdis) node [midway, right]  {$1$};
                \draw[-stealth] (80+0.1*\levdis, 0.1*\levdis) -- (40 + 0.1*\levdis, 0.9*\levdis) node [midway, right]  {$-\ratioobl{1}{8}$};
                
                \draw (40, \levdis) node [above] {$2k$};
                \draw (16, 0) node [below] {$k-1$};
                \draw (48, 0) node [below] {$k$};
                \draw (80, 0) node [below] {$k+1$};
                
            \end{tikzpicture}
        \end{center}
        \caption[]{\label{fig:prediction} Illustration of the action of the prediction operator taking $\predictionstencildepth = 1$. The cell average on the cell at level $\levelletter + 1$ is reconstructed by taking the values $\levelletter$ of the parent and his two neighbors with suitable weights.}
    \end{figure}

The opposite happens for the prediction operator (Figure \ref{fig:prediction}), taking information on a certain level $\levelletter$ and trying to recover an estimation of the values on a finer level $\levelletter + 1$. 
It seems that we have an infinity of possible choices and this is indeed the case.
However, we impose, following Cohen \emph{et al.} \cite{cohen2003}, some reasonable rigidity on the choice of the operator.

\begin{definition}[Prediction operator]
    The prediction operator $\predictionoperator : \reals^{1+w} \to \reals^2$, giving approximated values (denoted by a hat) of means at a fine level from data on a coarse level, that is
    \begin{equation*}
       \adaptiveroundbrackets{\predicted{\average{\distributionsadaptive}}{\populationletter}{\levelletter + 1, 2\indexletter + \delta}}_{\delta = 0, 1} = \predictionoperator \adaptiveroundbrackets{\adaptiveroundbrackets{\average{\distributionsadaptive}^{\populationletter}_{\levelletter, \pi}}_{\pi \in R(\levelletter, \indexletter)}},
    \end{equation*}
    for $\populationletter = 0, \dots, \velocitynumber - 1$, $\levelletter = \minlevel, \dots, \maxlevel - 1$ and $\indexletter = 0, \dots, N_{\levelletter}-1$, where
    \begin{itemize}
        \item The operator is local, namely the outcome depends on the value on $1+w$ cells at level $\levelletter$ with indices belonging to $R(\levelletter, \indexletter) $ geometrically close to $I_{\levelletter + 1, 2\indexletter + \delta}$ with $\delta  = 0, 1$.
        \item The operator is consistent with the projection operator, namely
        \begin{equation*}
            \projectionoperator \adaptiveroundbrackets{\adaptiveroundbrackets{\predicted{\average{\distributionsadaptive}}{\populationletter}{\levelletter + 1, 2\indexletter + \delta}}_{\delta = 0, 1}} = \average{\distributionsadaptive}_{\levelletter, \indexletter}^{\populationletter},
        \end{equation*}
        for $\populationletter = 0, \dots, \velocitynumber - 1$, $\levelletter = \minlevel, \dots, \maxlevel - 1$ and $\indexletter = 0, \dots, N_{\levelletter} - 1$.
    \end{itemize}
\end{definition}

\begin{remark}[Consequences of the definition]
    By the previous definition the parent belongs to the prediction stencil of its children (this is the $1$ in  $1+w$).
    Also observe that this definition does not impose to consider linear operators, even if it is the choice for this and many preceding works \cite{harten1994, harten1995, cohen2003, duarte2011,nguessan2019,nguessan2020}.
\end{remark}

In particular, let $\predictionstencildepth \in \mathbb{N}$ and consider for any $\populationletter = 0, \dots, \velocitynumber - 1$, $\levelletter = \minlevel, \dots, \maxlevel - 1$ and $\indexletter = 0, \dots, N_{\levelletter} - 1$
\begin{equation*}
    \predicted{\average{\distributionsadaptive}}{\populationletter}{\levelletter + 1, 2\indexletter + \delta} = \average{\distributionsadaptive}_{\levelletter, \indexletter}^{\populationletter} + (-1)^{\delta} \sum_{\alpha = 1}^{\gamma} c_{\alpha} \adaptiveroundbrackets{\average{\distributionsadaptive}_{\levelletter, \indexletter + \alpha}^{\populationletter} - \average{\distributionsadaptive}_{\levelletter, \indexletter - \alpha}^{\populationletter}}, \qquad \delta = 0, 1,
\end{equation*}
corresponding to the polynomial centered interpolations, which are exact for the averages of polynomials up to degree $2\predictionstencildepth$, being accurate at order $\mu \definitionequality 2 \predictionstencildepth + 1$.
Some coefficients are given by (see \cite{duarte2011,tenaud2011,nguessan2020} and references therein):
\begin{itemize}
    \item $\predictionstencildepth = 1$, with $c_1 = -\ratioobl{1}{8}$.
    \item $\predictionstencildepth = 2$, with $c_1 = -\ratioobl{22}{128}$ and $c_2 = \ratioobl{3}{128}$.
    \item $\predictionstencildepth = 3$, with $c_1 = -\ratioobl{201}{1024}$, $c_2 = \ratioobl{11}{256}$ and $c_3 = \ratioobl{5}{1024}$.
\end{itemize}

\subsection{Details and smoothness estimation}

Intuitively, the more the predicted value is far from the actual value on the considered cell, the more we can assume that the function locally lacks in smoothness due to the fact that it is far from behaving polynomially.
This is what is quantified by the notion of detail:

\begin{definition}[Detail]
    The details are defined as
    \begin{equation*}
        \average{\detailsadaptive}_{\levelletter, \indexletter}^{\populationletter} \definitionequality \average{\distributionsadaptive}_{\levelletter, \indexletter}^{\populationletter} - \predicted{\average{\distributionsadaptive}}{\populationletter}{\levelletter, \indexletter},
    \end{equation*}
    for $\populationletter = 0, \dots, \velocitynumber - 1$, $\levelletter = \minlevel + 1, \dots, \maxlevel$ \additionFree{and $\indexletter = 0, \dots, N_j-1$.}
\end{definition}

The details are redundant between siblings $I_{\levelletter + 1, 2\indexletter}$ and $I_{\levelletter + 1, 2\indexletter + 1}$ sharing the same parent $I_{\levelletter, \indexletter}$.
This is a trivial consequence of the consistency of $\predictionoperator$ and reads
\begin{equation}\label{eq:DetailRedundancy}
    \average{\detailsadaptive}_{\levelletter + 1, 2\indexletter}^{\populationletter} = - \average{\detailsadaptive}_{\levelletter + 1, 2\indexletter + 1}^{\populationletter},
\end{equation}
for $\populationletter = 0, \dots, \velocitynumber - 1$, $\levelletter = \minlevel, \dots, \maxlevel - 1$ and $k = 0, \dots, N_{j}-1$.
For this reason, we have to avoid the redundancy by considering only one detail between two siblings: we chose to keep only the one of the son with even indices $\average{\detailsadaptive}_{\levelletter + 1, 2\indexletter}^{\populationletter}$ \footnote{The opposite is perfectly fine.}.
Thus we introduce the sets of indices
\begin{align*}
    \setofindices_{\minlevel} &\definitionequality \adaptivecurlybrackets{(\minlevel, \indexletter) \quad : \quad \indexletter = 0, \dots, N_{\minlevel} - 1}, \\
    \setofindices_{\levelletter} &\definitionequality \adaptivecurlybrackets{(\levelletter, \indexletter) \quad : \quad k = 0, \dots, N_{\levelletter} - 1 \quad \text{and} \quad k ~ \text{even}}, \qquad \levelletter = \minlevel + 1, \dots, \maxlevel.
\end{align*}
Hence we have constructed the multiresolution $\multiresolutiontransform$ transform acting as follows
\begin{equation}\label{eq:actionMultiresolutionTransform}
    \vectorial{\average{\distributionsadaptive}}_{\maxlevel}^{\populationletter} \qquad \underset{\multiresolutiontransform^{-1}}{\overset{\multiresolutiontransform}{\myrightleftarrows{\rule{1cm}{0cm}}}} \qquad \adaptiveroundbrackets{\vectorial{\average{\distributionsadaptive}}_{\minlevel}^{\populationletter}, \vectorial{\average{\detailsadaptive}}_{\minlevel + 1}^{\populationletter}, \dots, \vectorial{\average{\detailsadaptive}}_{\maxlevel}^{\populationletter}},
\end{equation}
for every $\populationletter = 0, \dots, \velocitynumber - 1$, where
\begin{align*}
    \vectorial{\average{\distributionsadaptive}}_{\levelletter}^{\populationletter} &\definitionequality \adaptiveroundbrackets{\average{\distributionsadaptive}^{\populationletter}_{\levelletter, \indexletter}}_{\indexletter = 0, \dots, N_{\levelletter} - 1}, \qquad \levelletter = \minlevel, \dots, \maxlevel, \\
    \vectorial{\average{\detailsadaptive}}_{\levelletter}^{\populationletter} &\definitionequality \adaptiveroundbrackets{\average{\detailsadaptive}^{\populationletter}_{\levelletter, \indexletter}}_{(\levelletter, \indexletter) \in \setofindices_{\levelletter}}, \qquad \levelletter = \minlevel + 1, \dots, \maxlevel.
\end{align*}
One easily checks that each side of \eqref{eq:actionMultiresolutionTransform} contains the same number of elements, because we have eliminated the redundancy of the details.
\additionRefOne{The link between multiresolution and wavelets is the following, see \cite{cohen2003}: the averages are defined \emph{via} the ``dual-scaling function'' $\tilde{\varphi}_{\levelletter, \indexletter} \definitionequality \chi_{I_{\levelletter, \indexletter}}/|I_{\levelletter, \indexletter}|$, thus $\average{\distributionsadaptive}^{\populationletter}_{\levelletter, \indexletter} = \langle \distributionsadaptive^{\populationletter}, \tilde{\varphi}_{\levelletter, \indexletter}  \rangle$. Since we have selected linear prediction operators, they can be written as $\predicted{\average{\distributionsadaptive}}{\populationletter}{\levelletter, \indexletter} = \sum_{\tilde{\indexletter}} c_{\levelletter - 1, \tilde{\indexletter}}^{\levelletter, \indexletter} \average{\distributionsadaptive}^{\populationletter}_{\levelletter - 1, \tilde{\indexletter}}$ with suitable weights. Then
\begin{align*}
    \average{\detailsadaptive}_{\levelletter, \indexletter}^{\populationletter} \definitionequality \average{\distributionsadaptive}_{\levelletter, \indexletter}^{\populationletter} - \predicted{\average{\distributionsadaptive}}{\populationletter}{\levelletter, \indexletter} = \langle \distributionsadaptive^{\populationletter}, \tilde{\varphi}_{\levelletter, \indexletter}  \rangle - \sum_{\tilde{\indexletter}} c_{\levelletter - 1, \tilde{\indexletter}}^{\levelletter, \indexletter} \langle \distributionsadaptive^{\populationletter}, \tilde{\varphi}_{\levelletter - 1, \tilde{\indexletter}}  \rangle = \langle \distributionsadaptive^{\populationletter}, \tilde{\psi}_{\levelletter, \indexletter} \rangle,
\end{align*}
where we have introduced the dual wavelet $ \tilde{\psi}_{\levelletter, \indexletter} \definitionequality  \tilde{\varphi}_{\levelletter, \indexletter}  - \sum_{\tilde{\indexletter}} c_{\levelletter - 1, \tilde{\indexletter}}^{\levelletter, \indexletter}  \tilde{\varphi}_{\levelletter - 1, \tilde{\indexletter}} $. We shall indicate $\tilde{\Sigma}_{\levelletter, \indexletter} \definitionequality \text{supp}( \tilde{\psi}_{\levelletter, \indexletter})$ the support of the dual wavelet. The exactness of the prediction operators up to the order $2\predictionstencildepth$ can be interpreted as a vanishing property for the moments of the dual wavelet: if $\distributionsadaptive^{\populationletter} \in \Pi_{2\predictionstencildepth}$, the polynomials of degree at most $2\predictionstencildepth$, then $\average{\detailsadaptive}_{\levelletter, \indexletter}^{\populationletter} = \langle \distributionsadaptive^{\populationletter}, \tilde{\psi}_{\levelletter, \indexletter}\rangle = 0$, thus all the corresponding details are zero.
}
The details are a regularity indicator of the encoded function, as stated by the following Proposition
\begin{proposition}[Details decay]
    Consider a cell $I_{\levelletter, \indexletter}$ for $\levelletter = \minlevel + 1, \dots, \maxlevel$ and a population $\populationletter = 0, \dots, \velocitynumber - 1$ assuming that $\distributionsadaptive^{\populationletter} \in W_{\infty}^{\nu}(\tilde{\Sigma}_{\levelletter, \indexletter})$ for some $\nu \geq 0$, where
    \begin{equation*}
        W_{p}^{\nu}(I) \definitionequality \{\phi ~ : ~ \phi^{(\eta)} \in L^p(I), \quad 0 \leq \eta \leq \nu \}, \quad \lVert \phi \rVert_{W_{p}^{\nu}(I)} \definitionequality \lVert \phi \rVert_{L^p(I)} + |\phi|_{ W_{p}^{\nu}(I)},
    \end{equation*}
    where the semi-norm is $|\phi|_{ W_{p}^{\nu}(I)} \definitionequality \lVert \phi^{(\nu)}\rVert_{L^p(I)}$.
    Then, we have the following decay estimate for the details
    \begin{equation}\label{eq:DetailsDecayEstimate}
       |\average{\detailsadaptive}_{\levelletter, \indexletter}^{\populationletter}| \lesssim 2^{-\levelletter \min{(\nu, \mu)}} |\distributionsadaptive^{\populationletter}|_{W_{\infty}^{\min{(\nu, \mu)}}(\tilde{\Sigma}_{\levelletter, \indexletter})},
    \end{equation}
    \additionRefTwo{where the involved constants depend only on $\mu \definitionequality 2 \predictionstencildepth + 1$, where $\predictionstencildepth$ is the width of the prediction stencil.}
    \begin{proof}
        \additionRefOne{The proof is well-known in wavelet theory. It is provided in the Supplementary material for the interested reader.}
    \end{proof}
\end{proposition}

\begin{remark}
    Since the spaces $W_{\infty}^{\overline{\mu}}$ are algebr\ae{}, we can infer the regularity of the densities $\distributionsadaptive$ from the expected regularity of the moments \additionRefTwo{(obtained by the application of the matrix $\operatorial{M}$)}, in particular the conserved ones. \additionRefTwo{This is important because the conserved moments are eventually the quantities we are interested in and which have corresponding continuous equations under the form of conservation laws.}
\end{remark}
 This inequality\footnote{The interested reader can find a related numerical study in the Supplementary material.} states that the details become small when the function is locally smooth and also that they decrease with the level $\levelletter$ if functions are slightly more than just bounded.

\subsection{Tree structure and grading}

We introduce the set of all indices given by
\begin{equation*}
    \setofindices \definitionequality \bigcup_{\levelletter = \minlevel}^{\maxlevel} \setofindices_{\levelletter}.
\end{equation*}
In order to guarantee the feasibility of all the operations involved with the multiresolution and because  it naturally provides a multi-level covering of the domain $\Omega$, we want that our structure $\Lambda \subset \setofindices$ represents a graded tree.

\begin{definition}[Tree]
    Let $\Lambda \subset \setofindices$ be a set of indices. We say that $\Lambda$ represents a tree if
    \begin{enumerate}
        \item The coarsest level wholly belongs to the structure: $\setofindices_{\minlevel} \subset \Lambda$.
        \item There is no orphan cell: if $(\levelletter, \indexletter) \in \Lambda$, then $\adaptiveroundbrackets{\levelletter - 1, k/2} \in \Lambda$, for $\levelletter = \minlevel + 1, \dots, \maxlevel$\footnote{We are sure that $k/2$ is integer because we have only kept even cells}.
    \end{enumerate}
\end{definition}
Since we have discarded from $\setofindices$ the cells having redundant detail, given a tree $\Lambda \subset \setofindices$, we consider the complete tree $\physicaltree{\Lambda}$ obtained by adding the discarded cell with a siblings in $\Lambda$. With our choice, it is
\begin{equation*}
    \physicaltree{\Lambda} = \setofindices_{\minlevel} ~ \cup ~ \adaptivecurlybrackets{(\levelletter, \indexletter), (\levelletter, \indexletter + 1) ~ : ~ (\levelletter, \indexletter) \in \Lambda \quad \text{for} \quad \levelletter = \minlevel + 1, \dots, \maxlevel}.
\end{equation*}
Remark that $\Lambda \subsetneqq \physicaltree{\Lambda} \not\subset \setofindices$.
We also introduce the set of leaves $\leaves{\Lambda} \subset \Lambda \subset \setofindices$ of a tree $\Lambda$ which is the set of cells without child.
As usual, we introduce the set of complete leaves $\physicalleaves{\Lambda}$ which is given by
\begin{equation*}
    \physicalleaves{\Lambda} = \adaptivecurlybrackets{(\minlevel, k) ~ : ~ (\minlevel, \indexletter) \in \leaves{\Lambda}} ~\cup~ \adaptivecurlybrackets{(\levelletter, \indexletter), (\levelletter, \indexletter + 1) ~ : ~ (\levelletter, \indexletter) \in \leaves{\Lambda} \quad \text{and} \quad \levelletter > \minlevel}.
\end{equation*}
As observed by \cite{cohen2003}, the tree structure and the nesting allow us to conclude that $\physicalleaves{\Lambda}$ is a multi-level partition of the domain $\Omega$, \additionRefOne{therefore characterizing the locally refined mesh.}
Moreover $\leaves{\Lambda} \subsetneqq \physicalleaves{\Lambda} \not\subset \setofindices$.
We are ready to provide the definition of graded tree.
\begin{definition}[Graded tree]
    Let $\Lambda \subset \setofindices$ be a tree, then it is graded with respect to the prediction operator $\predictionoperator$ if the prediction stencil to predict over cells belonging to $\physicaltree{\Lambda} \smallsetminus \setofindices_{\minlevel}$ also belongs to $\physicaltree{\Lambda}$.
    With our prediction operator, this means that
    \begin{equation*}
        \text{If} \qquad (\levelletter, \indexletter) \in \physicaltree{\Lambda} \smallsetminus \setofindices_{\minlevel}, \quad \text{then} \quad (j - 1, \lfloor k/2 \rfloor + \delta) \in \physicaltree{\Lambda}, \qquad \delta = -\predictionstencilwidth, \dots, \predictionstencilwidth,
    \end{equation*}
    or equivalently, since we have removed redundant odd details
    \begin{equation*}
        \text{If} \qquad (\levelletter, \indexletter) \in \Lambda \smallsetminus \setofindices_{\minlevel}, \quad \text{then} \quad (j - 1, k/2  + \delta) \in \physicaltree{\Lambda}, \qquad \delta = -\predictionstencilwidth, \dots, \predictionstencilwidth.
    \end{equation*}
\end{definition}

Thus, given a tree $\Lambda \subset \setofindices$, we denote the operator yielding the smallest graded tree containing $\Lambda$ as $\gradingoperator(\Lambda)$.
The grading property is important because it guarantees that we can implement the isomorphism between
\begin{equation*}
    \adaptiveroundbrackets{\average{\distributionsadaptive}^{\populationletter}_{\levelletter, \indexletter}}_{(\levelletter, \indexletter) \in \physicalleaves{\Lambda}} \qquad \myrightleftarrows{\rule{1cm}{0cm}} \qquad \adaptiveroundbrackets{\vectorial{\average{\distributionsadaptive}}_{\minlevel}^{\populationletter}, \adaptiveroundbrackets{\average{\detailsadaptive}^{\populationletter}_{\levelletter, \indexletter}}_{(\levelletter, \indexletter) \in \Lambda \smallsetminus \setofindices_{\minlevel}}},
\end{equation*}
 for every $\populationletter = 0, \dots, \velocitynumber - 1$ in an efficient manner \cite{cohen2003}. This means that it is the same to know the means on the complete leaves $\physicalleaves{\Lambda}$ of a graded tree $\Lambda$ or knowing the averages on $\setofindices_{\minlevel} \subset \Lambda$ plus the details of $\Lambda \smallsetminus \setofindices_{\minlevel}$. In this work, we choose to store information on the complete leaves $\physicalleaves{\Lambda}$.

Let now $\Lambda \subset \setofindices$ be a graded tree and assume to know $(\average{\distributionsadaptive}^{\populationletter}_{\levelletter, \indexletter})_{(\levelletter, \indexletter) \in \physicalleaves{\Lambda}}$ or equivalently $(\vectorial{\average{\distributionsadaptive}}_{\minlevel}^{\populationletter}, (\average{\detailsadaptive}^{\populationletter}_{\levelletter, \indexletter})_{(\levelletter, \indexletter) \in \Lambda \smallsetminus \setofindices_{\minlevel}})$ for every $\populationletter = 0, \dots, \velocitynumber - 1$.
From this information, we can build the reconstruction on every cell at the finest level $\maxlevel$
\begin{equation*}
    \reconstructed{\vectorial{\average{\distributionsadaptive}}}{\populationletter}{\maxlevel} \definitionequality \adaptiveroundbrackets{\reconstructed{\average{\distributionsadaptive}}{\populationletter}{\maxlevel, \indexletter}}_{\indexletter = 0, \dots, N_{\maxlevel} - 1},
\end{equation*}
for $\populationletter = 0, \dots, \velocitynumber - 1$, where the double hat represents the reconstruction operator. With this operator, the information, stored on the complete leaves $\physicalleaves{\Lambda}$, is propagated from coarse (at the local level of resolution of $\physicalleaves{\Lambda}$) to the finest level $\maxlevel$ by means of level-by-level applications of the prediction operator  $\predictionoperator$.
\additionRefTwo{The reconstruction operator yields reconstructions of the lacking information on (possibly) virtual cells at the finest level (output) using the values stored on the complete leaves $\physicalleaves{\Lambda}$ at the local level of refinement (input).}

\subsection{Compressing information}
    
Take a graded tree $\Lambda \subset \setofindices$, then we consider the thresholding (or coarsening) operator given by
\begin{equation*}
    \thresholdoperator{\Lambda} \definitionequality \setofindices_{\minlevel} ~ \cup ~ \adaptivecurlybrackets{(\levelletter, \indexletter) \in \Lambda \smallsetminus \setofindices_{\minlevel} ~ : ~ \max_{\populationletter = 0, \dots, \velocitynumber - 1} |\average{\detailsadaptive}_{\levelletter, \indexletter}^{\populationletter}| \geq \thresholdletter_{\levelletter}} \subset \setofindices,
\end{equation*}
where the details concern $(\average{\distributionsadaptive}^{\populationletter}_{\levelletter, \indexletter})_{(\levelletter, \indexletter) \in \physicalleaves{\Lambda}}$ for $\populationletter = 0, \dots, \velocitynumber - 1$.
\additionRefTwo{This operator is constructed so that we end up with a compressed mesh which is the same for every field spanned by $\populationletter$ and is constructed by the most restrictive inequality on the details.}
It can be shown \cite{cohen2003, duarte2011} that
\begin{proposition}\label{prop:ErrorControl}
    Let $\thresholdletter > 0$ and consider a graded tree $\Lambda \subset \setofindices$ with data known on its complete leaves $\physicalleaves{\Lambda}$. 
    Consider the choice of level-wise thresholds
    \begin{equation*}
        \thresholdletter_{\levelletter} = 2^{\levelletter - \maxlevel} \epsilon, \qquad \levelletter = \minlevel + 1, \dots, \maxlevel,
    \end{equation*}
    and consider $p \in [0, \infty]$.
    Then there exists a constant $C_{\text{MR}} = C_{\text{MR}}(\gamma, p) > 0$ independent of $L$ such that
    \begin{equation*}
        \norm{\reconstructed{\vectorial{\average{\distributionsadaptive}}}{\populationletter}{\maxlevel} - \bm{A}_{\gradingoperator \circ \thresholdoperator{\Lambda}} \reconstructed{\vectorial{\average{\distributionsadaptive}}}{\populationletter}{\maxlevel}}_{\ell^p} \leq C_{\text{MR}} \epsilon,
    \end{equation*}
    for every $\populationletter = 0, \dots, \velocitynumber - 1$, where $\bm{A}_{\Lambda} \definitionequality \multiresolutiontransform^{-1} T_{\Lambda} \multiresolutiontransform$, where $T_{\Lambda}$ is the operator putting the details corresponding to indices which do not belong to $\Lambda$ at zero.
\end{proposition}

A similar estimate clearly holds when gathering all the populations spanning $\populationletter = 0, \dots, \velocitynumber - 1$.
It means that we can discard cells with small details still being able to reconstruct at the finest level $\maxlevel$ within a given precision controlled by $\thresholdletter$. We can then rely on this machinery in order to conduct a compression process and build a numerical strategy based on LBM schemes.
Observe that, from (43) in \cite{cohen2003}, the previous proposition also holds for $\vectorial{\average{\distributionsadaptive}}_{\maxlevel}^{\populationletter}$, thus also for $\reconstructed{\vectorial{\average{\distributionsadaptive}}}{\populationletter}{\maxlevel}$ as we stated.

\section{Adaptive MR-LBM scheme and error control}\label{sec:AdaptiveLBMMR}

So far, the procedure based on the multiresolution is static with respect to the evolution of time.
Since we want to utilize this strategy to build a reliable fully adaptive lattice Boltzmann solver for time dependent problems, it is of the foremost importance to define a way of evolving the compressed mesh so that it correctly represents the solution both at current time $t^n$ and at the successive time $t^{n+1}$, constructing it without \emph{a priori} knowing the new solution.

\subsection{Mesh adaptation strategy and time-stepping}

We are given an adaptive graded tree $\Lambda^n \subset \setofindices$ and a solution $(\average{\distributionsadaptive}_{\levelletter, \indexletter}^{\populationletter, n})_{(\levelletter, \indexletter) \in \physicalleaves{\Lambda^n}}$ for $\populationletter = 0, \dots, \velocitynumber - 1$ defined on the complete leaves of $\Lambda^n$ for the discrete time $t^n$.

\subsubsection{Mesh adaptation}

Starting from this level of information, which yields $(\vectorial{\average{\distributionsadaptive}}_{\minlevel}^{\populationletter, n}, (\average{\detailsadaptive}_{\levelletter, \indexletter}^{\populationletter, n})_{(\levelletter, \indexletter) \in \Lambda^n \smallsetminus \setofindices_{\minlevel}})$ using $\predictionoperator$ and $\projectionoperator$, we want to create a new mesh $\Lambda^{n+1}$ used to compute the new solution at time $t^{n+1}$ on $\physicalleaves{\Lambda^{n+1}}$.
The procedure can be schematized as follows
\begin{equation*}
    \Lambda^{n+1} \definitionequality \gradingoperator \circ \enlargingoperator \circ \thresholdoperator{\Lambda^n},
\end{equation*}
where the details used by $\enlargingoperator$ (still to be defined) and $\mathcal{T}_{\bm{\thresholdletter}}$ are those of the old solution, namely$(\average{\detailsadaptive}_{\levelletter, \indexletter}^{\populationletter, n})_{(\levelletter, \indexletter) \in \Lambda^n \smallsetminus \setofindices_{\minlevel}}$.
In the previous expression, we have:
\begin{itemize}
    \item $\mathcal{T}_{\bm{\thresholdletter}}$ is the thresholding operator we have previously defined. It can only merge fine cells on the tree to form coarser ones (coarsen).
    \item $\enlargingoperator$ is the enlargement operator. It breaks cells to form finer ones (refine) and is constructed to slightly enlarge the structure in order to accommodate the slowly evolving solution at the new time $t^{n+1}$.
    \item $\gradingoperator$ is the grading operator, which can also refine.
\end{itemize}
\additionRefTwo{We again observe that the details used to build $\mathcal{T}_{\bm{\thresholdletter}}$ and $\enlargingoperator$ are computed using the solution available at time $t^n$ on $\physicalleaves{\Lambda^n}$. The non-linear dependency of these operators on this solution is not written explicitly for the sake of keeping notations simple.}
Once we have $\Lambda^{n+1}$, we adapt the solution from $\physicalleaves{\Lambda^{n}}$ to $\physicalleaves{\Lambda^{n + 1}}$. When passing from $\Lambda^n$ to $\Lambda^{n+1}$, if cells are coarsened, we have to merge their data with the projection operator $\projectionoperator$. On the other hand, when finer cells are added by $\enlargingoperator$ or $\gradingoperator$, the missing information is reconstructed using the prediction operator $\predictionoperator$.
We are left with the old solution at time $t^n$ on the complete leaves of the new mesh $\physicalleaves{\Lambda^{n+1}}$: $(\average{\distributionsadaptive}_{\levelletter, \indexletter}^{\populationletter, n})_{(\levelletter, \indexletter)\in \physicalleaves{\Lambda^{n+1}}}$ for $\populationletter = 0, \dots, \velocitynumber - 1$.

\subsubsection{Time-stepping}

We denote the operator associated with the adaptive MR-LBM scheme by $\adaptiveschemeoperator^n$ (described in the sequel), with explicit dependence on the time $t^n$ since acting only on data defined on the time varying $\physicalleaves{\Lambda^{n+1}}$. It gives the approximate solution at time $t^{n+1}$ on the same hybrid grid.
This is
\begin{equation*}
    \adaptiveroundbrackets{\average{\distributionsadaptive}_{\levelletter, \indexletter}^{\populationletter, n+1}}_{\substack{(\levelletter, \indexletter)\in \physicalleaves{\Lambda^{n+1}} \\ \populationletter = 0, \dots, \velocitynumber - 1}} = \adaptiveschemeoperator^n \adaptiveroundbrackets{\average{\distributionsadaptive}_{\levelletter, \indexletter}^{\populationletter, n}}_{\substack{(\levelletter, \indexletter)\in \physicalleaves{\Lambda^{n+1}} \\ \populationletter = 0, \dots, \velocitynumber - 1}}.
\end{equation*}

\subsection{Construction of the enlargement operator $\mathcal{H}_{\epsilon}$}

We still have to define the enlargement operator $\mathcal{H}_{\epsilon}$, which is based on the following principles:

\begin{itemize}
    \item We must ensure that the propagation of information at finite speed \emph{via} the stencil of the lattice Boltzmann operator $\referenceschemeoperator$ (and thus also $\adaptiveschemeoperator^n$) is correctly handled.
    Thus, setting $\sigma = \max_{\populationletter =  0, \dots, \velocitynumber - 1}|\normalizedvelocityletter^{\populationletter}|$, we enforce that:
    \begin{equation*}
        \text{If} \quad (\levelletter, \indexletter) \in \physicaltree{\thresholdoperator{\Lambda^{n}}}, \quad \text{then} \quad (\levelletter, \indexletter + \delta) \in \physicaltree{\enlargingoperator \circ \thresholdoperator{\Lambda^{n}}}, \quad \delta = -\sigma, \dots, \sigma,
    \end{equation*}

    \item We must detect the shock formation possibly induced by the non linearity of the collisional part of $\referenceschemeoperator$ (or $\adaptiveschemeoperator^n$). Consider $\overline{\mu} \geq 0$ to be tuned, then
    \begin{gather}
        \text{If} \qquad (\levelletter, \indexletter) \in \thresholdoperator{\Lambda^n}, \quad \minlevel < \levelletter < \maxlevel, \quad \text{and} \quad  \max_{\populationletter = 0, \dots, \velocitynumber - 1} |\average{\detailsadaptive}_{\levelletter, \indexletter}^{\populationletter, n}| \geq 2^{\overline{\mu} + 1}\thresholdletter_{\levelletter}, \nonumber \\
        \text{then} \quad (\levelletter + 1, 2\indexletter + \delta) \in \physicaltree{\enlargingoperator \circ \thresholdoperator{\Lambda^{n}}}, \quad \delta = 0, 1, 2, 3. \label{eq:HartenLossOfRegularity}
    \end{gather}
    The rationale is the following: assume that the function $f^{\populationletter}(t^{n+1}, x)$ corresponding to $(\average{\distributionsadaptive}_{\levelletter, \indexletter}^{\populationletter, n+1})_{(\levelletter, \indexletter) \in \physicalleaves{\Lambda^{n+1}}}$ is such that $f^{\populationletter}(t^{n+1}, \cdot) \in W_{\infty}^{\nu}(\tilde{\Sigma}_{\levelletter, \indexletter})$ for some $\nu \geq 0$.
    Set $\overline{\mu} = \min{(\nu, \mu)}$.
    Since this solution is unknown at the stage at which we are utilizing $\enlargingoperator$, we assume that the solution varies slowly from $t^n$ to $t^{n+1}$, so that we claim 
    \begin{equation*}
        |\average{\detailsadaptive}_{\levelletter, \indexletter}^{\populationletter, n+1}| \approx |\average{\detailsadaptive}_{\levelletter, \indexletter}^{\populationletter, n}| \approx 2^{-\levelletter \overline{\mu}} |\distributionsadaptive^{\populationletter}(t^n, \cdot)|_{W_{\infty}^{\overline{\mu}}(\tilde{\Sigma}_{\levelletter, \indexletter})},
    \end{equation*}
    using \eqref{eq:DetailsDecayEstimate} and for the details which may not be available in the structure
    \begin{align*}
        |\average{\detailsadaptive}_{\levelletter + 1, 2\indexletter}^{\populationletter, n+1}| \approx |\average{\detailsadaptive}_{\levelletter + 1, 2\indexletter}^{\populationletter, n}| & \approx 2^{-(\levelletter + 1)\overline{\mu}} |\distributionsadaptive^{\populationletter}(t^n, \cdot)|_{W_{\infty}^{\overline{\mu}}(\tilde{\Sigma}_{\levelletter + 1, 2\indexletter})}, \\
        &\leq 2^{-(\levelletter + 1)\overline{\mu}} |\distributionsadaptive^{\populationletter}(t^n, \cdot)|_{W_{\infty}^{\overline{\mu}}(\tilde{\Sigma}_{\levelletter, \indexletter})},
    \end{align*}
    using the nesting of the lattices. 
    As a consequence, we have
    \begin{equation}\label{eq:DetailRatio}
        |\average{\detailsadaptive}_{\levelletter + 1, 2\indexletter}^{\populationletter, n+1}| \approx 2^{-\overline{\mu}}  |\average{\detailsadaptive}_{\levelletter, \indexletter}^{\populationletter, n}|.
    \end{equation}
    According to the analysis to construct the truncation operator $\thresholdoperator{\Lambda^{n+1}}$\footnote{$n+1$ because we are trying to anticipate the evolution of the solution.}, we would have kept $I_{\levelletter + 1, 2\indexletter}$ and $I_{\levelletter + 1, 2\indexletter + 1}$ if {$|\average{\detailsadaptive}_{\levelletter + 1, 2\indexletter}^{\populationletter, n+1}| \geq \thresholdletter_{\levelletter + 1} = 2 \thresholdletter_{\levelletter}$} and a priori also  $I_{\levelletter + 1, 2\indexletter + 2}$ and $I_{\levelletter + 1, 2\indexletter + 3}$, because their parent has a detail with the same absolute value of its sibling. 
    This comes back, using the previous estimate, at doing so whenever $|\average{\detailsadaptive}_{\levelletter, \indexletter}^{\populationletter, n}| \geq 2^{\overline{\mu} + 1} \thresholdletter_{\levelletter}$.
    Since the local regularity $\nu$ of the solution at each time step is unknown, $\overline{\mu} = \min(\nu, \mu)$ is a parameter of the simulation to be set.
\end{itemize}

\emph{Modulo} this operation on the mesh, which slightly enlarges the set of kept cells, we claim that the following heuristics, inspired by the works of Harten \cite{harten1994}, holds:
\begin{assumption}[Harten heuristics]
    The tree $\thresholdoperator{\Lambda^{n}}$ has been enlarged into a graded tree $\Lambda^{n+1} = \gradingoperator \circ \enlargingoperator \circ \thresholdoperator{\Lambda^n}$ such that for the chosen $p \in [1, \infty]$
    \begin{equation*}
        \norm{\reconstructed{\vectorial{\average{\distributionsadaptive}}}{n}{\maxlevel} - \bm{A}_{\Lambda^{n+1}} \reconstructed{\vectorial{\average{\distributionsadaptive}}}{n}{\maxlevel}}_{\ell^p} \leq C_{\text{MR}} \thresholdletter, \qquad \norm{\referenceschemeoperator \reconstructed{\vectorial{\average{\distributionsadaptive}}}{n}{\maxlevel} - \bm{A}_{\Lambda^{n+1}} (\referenceschemeoperator \reconstructed{\vectorial{\average{\distributionsadaptive}}}{n}{\maxlevel})}_{\ell^p} \leq C_{\text{MR}} \thresholdletter.
    \end{equation*}
\end{assumption}
The first assumption inequality is naturally fulfilled using the fact that $\thresholdoperator{\Lambda^n} \subset \Lambda^{n+1}$.
The second inequality is \additionRefOne{potentially} verified upon having enlarged the mesh using $\enlargingoperator$, which has been built considering how the scheme operator $\referenceschemeoperator$ acts on the solution.
It basically means that the mesh is suitable for well representing the solution obtained by applying the reference scheme to the adaptive solution at the previous time step $t^n$ reconstructed on the finest level.
\additionRefOne{Observe that we do not rigorously prove that this assumption holds for our refinement strategy $\mathcal{H}_{\epsilon}$. As for the Finite Volume scheme, the Harten's approach to construct $\mathcal{H}_{\epsilon}$ has never proved to satisfy the assumption but is widely used in practice. The only achievement in terms of reliability condition has been obtained in \cite{cohen2003} for Finite Volume with scalar conservation laws, with a quite sophisticated refinement strategy.}

\subsection{Construction of the adaptive MR-LBM scheme}

We now present how to construct the adaptive MR-LBM  scheme $\adaptiveschemeoperator^n$ from the reference scheme $\referenceschemeoperator$.

\subsubsection{Collision}
In this part, the change of variable \emph{via} $\operatorial{M}$ is understood.
Consider a complete leaf $(\levelletter, \indexletter) \in \physicalleaves{\Lambda^{n+1}}$.
Consider all the cells of $\mathcal{L}_{\maxlevel}$ inside $I_{\levelletter, \indexletter}$ with indices $(\maxlevel, k2^{\maxlevel - \levelletter} + \delta)$ with $\delta = 0, \dots, 2^{\maxlevel - \levelletter} - 1$. We reconstruct the data on them and we project back on $I_{\levelletter, \indexletter}$ belonging to the complete leaves $\physicalleaves{\Lambda^{n+1}}$. This is
\begin{align}
        \average{\momentumadaptive}_{\levelletter, \indexletter}^{\populationletter, n\star} &= \average{\momentumadaptive}_{\levelletter, \indexletter}^{\populationletter, n}, \qquad \populationletter = 0, \dots, \velocitynumber_{\text{cons}} - 1, \nonumber \\ 
       \average{\momentumadaptive}_{\levelletter, \indexletter}^{\populationletter, n\star} &= (1-s^h)  \average{\momentumadaptive}_{\levelletter, \indexletter}^{\populationletter, n} + \ratio{s^h}{2^{\maxlevel - \levelletter}} \sum_{\delta = 0}^{2^{\maxlevel - \levelletter} - 1} M^{\populationletter, \text{eq}} \adaptiveroundbrackets{\reconstructed{\average{\momentumadaptive}}{0, n}{\maxlevel, \indexletter 2^{\maxlevel - \levelletter} + \delta}, \dots, \reconstructed{\average{\momentumadaptive}}{\velocitynumber_{\text{cons}} - 1, n}{\maxlevel, \indexletter 2^{\maxlevel - \levelletter} + \delta}}, \label{eq:reconstructedCollision}
\end{align}
for $\populationletter = \velocitynumber_{\text{cons}}, \dots, \velocitynumber - 1$.
\begin{remark}\label{rem:TooCherCollision}
As observed by \cite{hovhannisyan2010} for the source terms of Finite Volume schemes, this strategy can be computationally expensive and is mostly of theoretical interest.
We shall discuss this fact and introduce an alternative approach in the sequel.
\end{remark}
The first term on the right is taken on the complete leaf because it is linear, thus the reconstructed values simplify when taking the projection operator.

\subsubsection{Stream}

For the sake of notation, let us introduce the sign of each velocity given by $\sigma^{\populationletter} \definitionequality w^{\populationletter}/|w^{\populationletter}| \in \{-1, 0, 1 \}$ for $\populationletter = 0, \dots, \velocitynumber - 1$ fixed in the sequel. Consider a complete leaf $(\levelletter, \indexletter) \in \physicalleaves{\Lambda^{n+1}}$.
Consider all the cells of $\mathcal{L}_{\maxlevel}$ inside $I_{\levelletter, \indexletter}$ with indices $(\maxlevel, k2^{\maxlevel - \levelletter} + \delta)$ with $\delta = 0, \dots, 2^{\maxlevel - \levelletter} - 1$.
Perform the advection on them
\begin{equation*}
    \average{\distributionsadaptive}_{\maxlevel, k2^{\maxlevel - \levelletter} + \delta}^{\populationletter, n+1} = \average{\distributionsadaptive}_{\maxlevel, k2^{\maxlevel - \levelletter} + \delta - w^{\populationletter}}^{\populationletter, n\star}, \quad \text{for} \quad \delta = 0, \dots, 2^{\maxlevel - \levelletter} - 1.
\end{equation*}
The problem is that the data on the right hand side are usually unavailable since since their cell does not belong to $\physicalleaves{\Lambda^{n+1}}$.
Despite this, we can reconstruct, yielding
\begin{equation*}
    \average{\distributionsadaptive}_{\maxlevel, k2^{\maxlevel - \levelletter} + \delta}^{\populationletter, n+1} \approx
    \reconstructed{\average{\distributionsadaptive}}{\populationletter, n\star}{\maxlevel, k2^{\maxlevel - \levelletter} + \delta - w^{\populationletter}}, \quad \text{for} \quad \delta = 0, \dots, 2^{\maxlevel - \levelletter} - 1.
\end{equation*}
We want to update the solution on $\physicalleaves{\Lambda^{n+1}}$, that is why we project using  the projection operator $\projectionoperator$ $\maxlevel - \levelletter$ times
\begin{equation}\label{eq:StreamUncompacted}
    \average{\distributionsadaptive}_{\levelletter, k}^{\populationletter, n+1} \approx \ratio{1}{2^{\maxlevel - \levelletter}} \sum_{\delta = 0}^{2^{\maxlevel - \levelletter} - 1} \reconstructed{\average{\distributionsadaptive}}{\populationletter, n\star}{\maxlevel, k2^{\maxlevel - \levelletter} + \delta - w^{\populationletter}}.
\end{equation}
Indeed, only the terms referring to the virtual cells close to the boundary of the leaf are actually needed.
This can be seen in the following way. 

If $\sigma^{\populationletter} = 0$, then by the consistency of $\predictionoperator$ with $\projectionoperator$, \eqref{eq:StreamUncompacted} becomes
    \begin{equation*}
        \ratio{1}{2^{\maxlevel - \levelletter}} \sum_{\delta = 0}^{2^{\maxlevel - \levelletter} - 1} \reconstructed{\average{\distributionsadaptive}}{\populationletter, n\star}{\maxlevel, k2^{\maxlevel - \levelletter} + \delta - w^{\populationletter}}=  \average{\distributionsadaptive}_{\levelletter, k}^{\populationletter, n\star}.
    \end{equation*}

If $\sigma^{\populationletter} = 1 > 0$, then using consistency once more
    \begin{equation*}
        \ratio{1}{2^{\maxlevel - \levelletter}}  \sum_{\delta = 0}^{2^{\maxlevel - \levelletter} - 1} \reconstructed{\average{\distributionsadaptive}}{\populationletter, n\star}{\maxlevel, k2^{\maxlevel - \levelletter} + \delta - w^{\populationletter}} =  \average{\distributionsadaptive}_{\levelletter, k}^{\populationletter, n\star}+ \ratio{1}{2^{\maxlevel - \levelletter}} \sum_{\delta = 1}^{w^{\populationletter}} \left ( \reconstructed{\average{\distributionsadaptive}}{\populationletter, n\star}{\maxlevel, k2^{\maxlevel - \levelletter} - \delta} - \reconstructed{\average{\distributionsadaptive}}{\populationletter, n\star}{\maxlevel, (k+1)2^{\maxlevel - \levelletter} - \delta} \right ).
    \end{equation*}

If $\sigma^{\populationletter} = -1 < 0$, then by consistency
    \begin{equation*}
        \ratio{1}{2^{\maxlevel - \levelletter}}  \sum_{\delta = 0}^{2^{\maxlevel - \levelletter} - 1} \reconstructed{\average{\distributionsadaptive}}{\populationletter, n\star}{\maxlevel, k2^{\maxlevel - \levelletter} + \delta - w^{\populationletter}} =  \average{\distributionsadaptive}_{\levelletter, k}^{\populationletter, n\star} - \ratio{1}{2^{\maxlevel - \levelletter}} \sum_{\delta = 1}^{-w^{\populationletter}} \left ( \reconstructed{\average{\distributionsadaptive}}{\populationletter, n\star}{\maxlevel, k2^{\maxlevel - \levelletter} + \delta - 1} - \reconstructed{\average{\distributionsadaptive}}{\populationletter, n\star}{\maxlevel, (k+1)2^{\maxlevel - \levelletter} + \delta - 1} \right ).
    \end{equation*}
To unify all this formul\ae, we introduce $\eta(\populationletter, \delta) \in  \mathbb{Z}$ for $\delta = 1, \dots, |w^{\populationletter}|$ defined by
\begin{equation*}
    \eta(\populationletter, \delta) \definitionequality \left ( \frac{1}{2} - \delta \right ) \sigma^{\populationletter} - \frac{1}{2} = 
    \begin{cases}
        -\delta, \qquad &\text{if} \quad \sigma^{\populationletter} = 1 > 0, \\
        \delta - 1, \qquad &\text{if} \quad \sigma^{\populationletter} = -1 < 0. \\
    \end{cases}
\end{equation*}
Thus we obtain the more compact expression for the stream phase
\begin{equation}\label{eq:AdaptiveStreamPhase}
    \average{\distributionsadaptive}_{\levelletter, k}^{\populationletter, n+1} \approx \average{\distributionsadaptive}_{\levelletter, k}^{\populationletter, n\star} + \ratio{\sigma^{\populationletter}}{2^{\maxlevel - \levelletter}} \sum_{\delta = 1}^{|w^{\populationletter}|} \adaptiveroundbrackets{\reconstructed{\average{\distributionsadaptive}}{\populationletter, n\star}{\maxlevel, k2^{\maxlevel - \levelletter}+\eta(\populationletter, \delta)} - \reconstructed{\average{\distributionsadaptive}}{\populationletter, n\star}{\maxlevel, (k+1)2^{\maxlevel - \levelletter}+\eta(\populationletter, \delta)}}.
\end{equation}

\begin{remark}\label{rem:TooCherStream}
Since the reconstruction operator $\reallywidedoublehat{\hphantom{f}}$ which utilizes $\predictionoperator$ until reaching available values on $\physicalleaves{\Lambda^{n+1}}$ does not depend on the details (they are not available) one might use \additionFree{a} cheaper interpolation to perform this operation, as hinted by \cite{cohen2003}. Such an approach is used in many works, but at the cost of the error control provided by the MR machinery.
\end{remark}

\subsection{Error analysis}

The major interest of adaptive meshes generated by multiresolution is that we can recover a precise error control on the perturbation (or additional) error when solving PDEs on them.
Fixing a given $\ell^p$ norm for $p \in [1, \infty]$, we want to control the additional error $\lVert \vectorial{\average{F}}^n - \reconstructed{\vectorial{\average{f}}}{n}{\additionRefTwoSecondRound{\maxlevel}} \rVert_{\ell^p}$ where $\vectorial{\average{F}}^n$ is the solution of the reference scheme given by $\vectorial{\average{ F}}^{n+1} = \referenceschemeoperator \vectorial{\average{F}}^n$ and wholly defined on the finest level $\maxlevel$ and computations start from the same initial datum on the finest grid, that is $\Lambda^0 = \setofindices$.
In the following analysis, the assumptions are the following
\begin{itemize}
    \item \textbf{H1 - Harten heuristics}. At each step, the tree $\thresholdoperator{\Lambda^n}$ has been enlarged into a graded tree $\Lambda^{n+1}$ so that 
    \begin{equation*}
        \norm{\reconstructed{\vectorial{\average{\distributionsadaptive}}}{n}{\maxlevel} - \bm{A}_{\Lambda^{n+1}} \reconstructed{\vectorial{\average{\distributionsadaptive}}}{n}{\maxlevel}}_{\ell^p} \leq C_{\text{MR}} \thresholdletter, \qquad \norm{\referenceschemeoperator \reconstructed{\vectorial{\average{\distributionsadaptive}}}{n}{\maxlevel} - \bm{A}_{\Lambda^{n+1}} (\referenceschemeoperator \reconstructed{\vectorial{\average{\distributionsadaptive}}}{n}{\maxlevel})}_{\ell^p} \leq C_{\text{MR}} \thresholdletter.
    \end{equation*}
    \item \textbf{H2 - Continuity of $\referenceschemeoperator$}. There exists a constant $C_{L} = 1 + \tilde{C}_L$ with $\tilde{C}_L \geq 0$ such that
    \begin{equation*}
        \norm{\referenceschemeoperator \vectorial{\average{U}} - \referenceschemeoperator \vectorial{\average{V}}}_{\ell^p} \leq C_L \norm{\vectorial{\average{U}} - \vectorial{\average{V}}}_{\ell^p}, \qquad \forall \vectorial{\average{U}}, \vectorial{\average{V}} \in \reals^{\velocitynumber N_{\maxlevel}}.
    \end{equation*}
\end{itemize}
\begin{remark}
        The following procedure can be easily adapted to the context where the continuity of the scheme is measured using a $\ell^2$-weighted norm as by \cite{junk2009}. It is sufficient to consider $p = 2$ and to observe that the corresponding norm (measuring the properties pertaining to the multiresolution) can be bounded by the $\ell^2$-weighted norm.
    \end{remark}
    
Thus we prove, replicating the path of \cite{cohen2003}, the following statement which gives a control on the error introduced by the MR-LBM adaptive scheme\strikeformula{:}
Remark that our formulation of the Harten heuristics is slightly different from the one in \cite{cohen2003} because of the different order of the operations at each time step of the algorithm, see Section 3.3 in \cite{cohen2003}. However, this does not make any difference, except when dealing with the initial datum, because the order of the operations when time steps are concatenated is the same. In our work, we do not aim at providing a construction of the enlargement operator $\enlargingoperator$ that rigorously guarantees the Harten heuristics, contrarily to \cite{cohen2003}. We just state and assume the heuristics which shall eventually turn out to be fulfilled in the examples of section \ref{sec:Verifications}.

\begin{proposition}[Additional error estimate]\label{proposition:ErrorEstimate}
    Under the Assumptions (H1) and (H2), the additional error satisfies the following upper bounds
        \begin{equation*}
        \norm{\vectorial{\average{F}}^n - \reconstructed{\vectorial{\average{f}}}{n}{\additionRefTwoSecondRound{\maxlevel}}}_{\ell^p} \leq C_{\text{MR}} ~  \thresholdletter \times 
        \begin{cases}
            n+1, \qquad &\text{if} \quad \tilde{C}_L = 0, \\
            1 + \ratio{e^{\tilde{C}_L n } - 1}{\tilde{C}_L}, \qquad &\text{if} \quad \tilde{C}_L > 0. \\
        \end{cases}
    \end{equation*}
    \begin{proof}
        Start by observing that as stated in the proof of Proposition 4.2 in \cite{cohen2003} or (3.117) in \cite{duarte2013}, since we reconstruct at the finest level both for the collision and the stream phase
        \begin{equation}\label{eq:tmp0}
            \reconstructed{\vectorial{\average{\distributionsadaptive}}}{n+1}{\additionRefTwoSecondRound{\maxlevel}} = (\bm{A}_{\Lambda^{n+1}} \circ \referenceschemeoperator) \reconstructed{\vectorial{\average{\distributionsadaptive}}}{n}{\additionRefTwoSecondRound{\maxlevel}},
        \end{equation}
        where $\reconstructed{\vectorial{\average{\distributionsadaptive}}}{n}{\additionRefTwoSecondRound{\maxlevel}}$ is reconstructed from data already adapted on $\Lambda^{n+1}$.
        Hence by the triangle inequality
        \begin{align*}
            \lVert \vectorial{\average{F}}^{n} - \reconstructed{\vectorial{\average{f}}}{n}{\additionRefTwoSecondRound{\maxlevel}} \rVert_{\ell^p} &\leq \lVert \referenceschemeoperator\vectorial{\average{F}}^{n - 1} - \referenceschemeoperator \reconstructed{\vectorial{\average{f}}}{n-1}{\additionRefTwoSecondRound{\maxlevel}} \rVert_{\ell^p} +  \lVert \referenceschemeoperator \reconstructed{\vectorial{\average{f}}}{n-1}{\additionRefTwoSecondRound{\maxlevel}} - \reconstructed{\vectorial{\average{f}}}{n}{\additionRefTwoSecondRound{\maxlevel}} \rVert_{\ell^p},\\
            &\leq (1+\tilde{C}_L) \lVert\vectorial{\average{F}}^{n - 1} - \reconstructed{\vectorial{\average{f}}}{n-1}{\additionRefTwoSecondRound{\maxlevel}} \rVert_{\ell^p} + \lVert \referenceschemeoperator \reconstructed{\vectorial{\average{f}}}{n-1}{\additionRefTwoSecondRound{\maxlevel}} - (\bm{A}_{\Lambda^{n}} \circ \referenceschemeoperator) \reconstructed{\vectorial{\average{\distributionsadaptive}}}{n - 1}{\additionRefTwoSecondRound{\maxlevel}} \rVert_{\ell^p},\\
            &\leq (1+\tilde{C}_L) \lVert\vectorial{\average{F}}^{n - 1} - \reconstructed{\vectorial{\average{f}}}{n-1}{\additionRefTwoSecondRound{\maxlevel}} \rVert_{\ell^p} + C_{\text{MR}} \thresholdletter,
        \end{align*}
        employing in this order Assumption (H2), \eqref{eq:tmp0} and Assumption (H1).
        We have to distinguish two cases and apply the inequality recursively
        \begin{itemize}
            \item $\tilde{C}_L = 0$, thus $\lVert \vectorial{\average{F}}^{n} - \reconstructed{\vectorial{\average{f}}}{n}{\additionRefTwoSecondRound{\maxlevel}} \rVert_{\ell^p} \leq \lVert \vectorial{\average{F}}^{n-1} - \reconstructed{\vectorial{\average{f}}}{n-1}{\additionRefTwoSecondRound{\maxlevel}} \rVert_{\ell^p} + C_{\text{MR}}\thresholdletter \leq \dots \leq C_{\text{MR}} (n+1)\epsilon$. Observe that $n+1$ comes from the fact that $\lVert \vectorial{\average{F}}^{0} - \reconstructed{\vectorial{\average{f}}}{0}{\additionRefTwoSecondRound{\maxlevel}} \rVert_{\ell^p} \neq 0$, but we only have $\lVert \vectorial{\average{F}}^{0} - \reconstructed{\vectorial{\average{f}}}{0}{\additionRefTwoSecondRound{\maxlevel}} \rVert_{\ell^p} \leq C_{\text{MR}} \thresholdletter$ because of Proposition \ref{prop:ErrorControl}.
            \item $\tilde{C}_L > 0$. We obtain, using that $(1+\tilde{C})^n \leq e^{\tilde{C}n}$ if $\tilde{C} > 0$
            \begin{align*}
                \lVert \vectorial{\average{F}}^{n} - \reconstructed{\vectorial{\average{f}}}{n}{\additionRefTwoSecondRound{\maxlevel}} \rVert_{\ell^p} &\leq (1+\tilde{C}_L) \lVert\vectorial{\average{F}}^{n - 1} - \reconstructed{\vectorial{\average{f}}}{n-1}{\additionRefTwoSecondRound{\maxlevel}} \rVert_{\ell^p} + C_{\text{MR}} \thresholdletter \leq \dots \\
                &\leq C_{\text{MR}} \thresholdletter \sum_{i = 0}^{n-1} (1+\tilde{C}_L)^i +  C_{\text{MR}} \thresholdletter = C_{\text{MR}} \adaptiveroundbrackets{1 + \ratio{(1+\tilde{C}_L)^n - 1}{\tilde{C}_L}} \thresholdletter \\
                &\leq C_{\text{MR}} \adaptiveroundbrackets{1 + \ratio{e^{\tilde{C}_L n } - 1}{\tilde{C}_L}} \thresholdletter.
            \end{align*}
        \end{itemize}
    \end{proof}
\end{proposition}
Therefore, regardless of continuity constant of the reference scheme, the additional error is bounded linearly with $\thresholdletter$.
According to the value of constant, we can prove that it accumulates either at most linearly in time or exponentially.
\additionRefTwo{It is in general difficult to link the relaxation parameters with the constant and,  according to our experience, experiments frequently show a linear behavior even when we expect an exponential one, thus the bound is not sharp.}

\subsection{Conclusion, discussion and implementation}

    \begin{figure}
        \begin{center}
            \begin{tikzpicture}[x=0.04cm, y=0.04cm]
                \newcount\levdis
                        \levdis = 12

                        \draw [Bar-Bar, black, thick] (0,0) -- (16,0);
                        \draw [-Bar, black, thick] (16,0) -- (32,0);
                        \draw [-Bar, gray!50, dashed] (32,0) -- (48,0);
                        \draw [-Bar, gray!50, dashed] (48,0) -- (64,0);
                        \draw [-Bar, gray!50, dashed] (64,0) -- (80,0);
                        \draw [-, gray!50, dashed] (80,0) -- (96,0);
                        \draw [Bar-Bar, black, thick] (96,0) -- (112,0);
                        \draw [-Bar, black, thick] (112,0) -- (128,0);
                        \draw (0,0) node [left] {\footnotesize{$\maxlevel - 4$}};

                        \draw [Bar-Bar, black, thick] (32, \levdis) -- (40, \levdis);
                        \draw [-Bar, gray!50, dashed] (40, \levdis) -- (48, \levdis);
                        \draw [-Bar, gray!50, dashed] (48, \levdis) -- (56, \levdis);
                        \draw [-Bar, gray!50, dashed] (56, \levdis) -- (64, \levdis);
                        \draw [-Bar, gray!50, dashed] (64, \levdis) -- (72, \levdis);
                        \draw [-Bar, gray!50, dashed] (72, \levdis) -- (80, \levdis);
                        \draw [-, gray!50, dashed] (80, \levdis) -- (88, \levdis);
                        \draw [Bar-Bar, black, thick] (88, \levdis) -- (96, \levdis);
                        \draw (0,\levdis) node [left] {\footnotesize{$\maxlevel - 3$}};

                        \draw [Bar-Bar, black, thick] (40, 2*\levdis) -- (44, 2*\levdis);
                        \draw [-Bar, black, thick] (44, 2*\levdis) -- (48, 2*\levdis);
                        \draw [-Bar, black, thick] (48, 2*\levdis) -- (52, 2*\levdis);
                        \draw [-Bar, gray!50, dashed] (52, 2*\levdis) -- (56, 2*\levdis);
                        \draw [-Bar, gray!50, dashed] (56, 2*\levdis) -- (60, 2*\levdis);
                        \draw [-Bar, gray!50, dashed] (60, 2*\levdis) -- (64, 2*\levdis);
                        \draw [-Bar, gray!50, dashed] (64, 2*\levdis) -- (68, 2*\levdis);
                        \draw [-Bar, gray!50, dashed] (68, 2*\levdis) -- (72, 2*\levdis);
                        \draw [-, gray!50, dashed] (72, 2*\levdis) -- (76, 2*\levdis);
                        \draw [Bar-Bar, black, thick] (76, 2*\levdis) -- (80, 2*\levdis);
                        \draw [-Bar, black, thick] (80, 2*\levdis) -- (84, 2*\levdis);
                        \draw [-Bar, black, thick] (84, 2*\levdis) -- (88, 2*\levdis);
                        \draw (0,2*\levdis) node [left] {\footnotesize{$\maxlevel - 2$}};
                        
                        \draw [densely dotted] (82,2*\levdis) circle (4) ;

                        \draw[->,>=latex] (82 + 3, 2*\levdis + 3) to[bend left=45]  (150, 60);
                        \draw[] (110, 55) node [right] {\footnotesize{$\forall$ leaf}};

                        \draw [Bar-Bar, black, very thick] (150, 5*\levdis) -- (214, 5*\levdis);
                        \path[>=latex, thick, densely dotted] (182, 64) edge [out=120, in=45, minimum size=16mm, loop] (182, 64);
                        
                        \draw[decoration={brace,raise=5pt},decorate]
                            (255,70) -- node[right=6pt] {\footnotesize{\textbf{Collision}}} (255, 50);
                        
                        \draw[->,>=latex] (82 + 3, 2*\levdis + 3) to[bend left=45]  (150, 40);
                        \draw [Bar-Bar, black, very thick] (150, 40) -- (214, 40);
                        
                        \draw[->,>=latex] (214, 36) to[bend left=45]  (214, 24);
                        \draw[] (218, 30) node [right] {\footnotesize{Split}};

                        \draw [Bar-, black, dashed] (134, 20) -- (150, 20);
                        \draw [Bar-Bar, black, very thick] (150, 20) -- (166, 20);
                        \draw [-Bar, black, very thick] (166, 20) -- (182, 20);
                        \draw [-Bar, black, very thick] (182, 20) -- (198, 20);
                        \draw [-Bar, black, very thick] (198, 20) -- (214, 20);
                        \draw [-Bar, black, dashed] (214, 20) -- (230, 20);
                        
                        \draw[decoration={brace,raise=5pt},decorate]
                            (255,45) -- node[right=6pt] {\footnotesize{\textbf{Stream}}} (255, -5);

                        \draw[->,>=latex, thick, densely dotted] (142, 20) to[bend left=75]  (158, 20);
                        \draw[->,>=latex, thick, densely dotted] (158, 20) to[bend left=75]  (174, 20);
                        \draw[->,>=latex, thick, densely dotted] (174, 20) to[bend left=75]  (190, 20);
                        \draw[->,>=latex, thick, densely dotted] (190, 20) to[bend left=75]  (206, 20);
                        
                        \draw[->,>=latex, thick, densely dotted] (222, 20) to[bend left=75]  (206, 20);
                        \draw[->,>=latex, thick, densely dotted] (206, 20) to[bend left=75]  (190, 20);
                        \draw[->,>=latex, thick, densely dotted] (190, 20) to[bend left=75]  (174, 20);
                        \draw[->,>=latex, thick, densely dotted] (174, 20) to[bend left=75]  (158, 20);
                        \draw[->,>=latex, thick, densely dotted] (158, 20) to[bend left=75]  (142, 20);
                        
                        \draw [Bar-Bar, black, very thick] (150, 0) -- (214, 0);
                        \draw[->,>=latex] (214, 16) to[bend left=45]  (214, 4);
                        \draw[] (218, 10) node [right] {\footnotesize{Projection}};

                        \draw [Bar-Bar, black, thick] (52, 3*\levdis) -- (54, 3*\levdis);
                        \draw [-Bar, black, thick] (54, 3*\levdis) -- (56, 3*\levdis);
                        \draw [-Bar, gray!50, dashed] (56, 3*\levdis) -- (58, 3*\levdis);
                        \draw [-, gray!50, dashed] (58, 3*\levdis) -- (60, 3*\levdis);
                        \draw [Bar-Bar, black, thick] (60, 3*\levdis) -- (62, 3*\levdis);
                        \draw [-Bar, gray!50, dashed] (62, 3*\levdis) -- (64, 3*\levdis);
                        \draw [-, gray!50, dashed] (64, 3*\levdis) -- (66, 3*\levdis);
                        \draw [Bar-Bar, black, thick] (66, 3*\levdis) -- (68, 3*\levdis);
                        \draw [-Bar, gray!50, dashed] (68, 3*\levdis) -- (70, 3*\levdis);
                        \draw [-, gray!50, dashed] (70, 3*\levdis) -- (72, 3*\levdis);
                        \draw [Bar-Bar, black, thick] (72, 3*\levdis) -- (74, 3*\levdis);
                        \draw [-Bar, black, thick] (74, 3*\levdis) -- (76, 3*\levdis);
                        \draw (0,3*\levdis) node [left] {\footnotesize{$\maxlevel - 1$}};

                        \draw [Bar-Bar, black, thick] (56, 4*\levdis) -- (57, 4*\levdis);
                        \draw [-Bar, black, thick] (57, 4*\levdis) -- (58, 4*\levdis);
                        \draw [-Bar, black, thick] (58, 4*\levdis) -- (59, 4*\levdis);
                        \draw [-Bar, black, thick] (59, 4*\levdis) -- (60, 4*\levdis);

                        \draw [Bar-Bar, black, thick] (62, 4*\levdis) -- (63, 4*\levdis);
                        \draw [-Bar, black, thick] (63, 4*\levdis) -- (64, 4*\levdis);
                        \draw [-Bar, black, thick] (64, 4*\levdis) -- (65, 4*\levdis);
                        \draw [-Bar, black, thick] (65, 4*\levdis) -- (66, 4*\levdis);

                        \draw [Bar-Bar, black, thick] (68, 4*\levdis) -- (69, 4*\levdis);
                        \draw [-Bar, black, thick] (69, 4*\levdis) -- (70, 4*\levdis);
                        \draw [-Bar, black, thick] (70, 4*\levdis) -- (71, 4*\levdis);
                        \draw [-Bar, black, thick] (71, 4*\levdis) -- (72, 4*\levdis);
                        \draw (0,4*\levdis) node [left] {\footnotesize{$\maxlevel$}};

                \end{tikzpicture}
        \end{center}
        \caption[]{\label{fig:standcoll} Schematic illustration of the basic features of the adaptive numerical scheme as it is currently used in the paper, meaning with the ``leaves collision'' \eqref{eq:LeavesCollision}.}
    \end{figure}
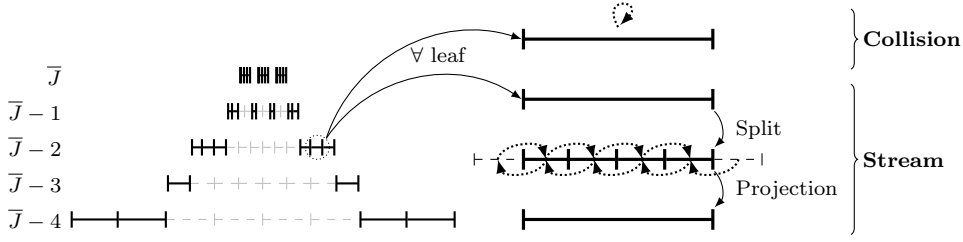


As we observed with Remark \ref{rem:TooCherCollision}, the collision given by \eqref{eq:reconstructedCollision} (called ``reconstructed collision'') is used in the theoretical analysis but remains limiting in practice, especially in the multidimensional context \cite{bellotti2021}.
Therefore, we propose the so-called ``leaves collision'' (see Figure \ref{fig:standcoll}), using data available on the complete leaves $\physicalleaves{\Lambda^{n+1}}$, \additionRefOne{which form the locally refined mesh}. This reads, for $(\levelletter, \indexletter) \in \physicalleaves{\Lambda^{n+1}}$
\begin{align}
        \average{\momentumadaptive}_{\levelletter, \indexletter}^{\populationletter, n\star} = \average{\momentumadaptive}_{\levelletter, \indexletter}^{\populationletter, n}, \quad &\populationletter = 0, \dots, \velocitynumber_{\text{cons}} - 1, \nonumber \\ 
        \average{\momentumadaptive}_{\levelletter, \indexletter}^{\populationletter, n\star} = (1-s^h)  \average{\momentumadaptive}_{\levelletter, \indexletter}^{\populationletter, n}  + s^h M^{\populationletter, \text{eq}} \adaptiveroundbrackets{\average{\momentumadaptive}_{\levelletter, \indexletter}^{0, n}, \dots, \average{\momentumadaptive}_{\levelletter, \indexletter}^{\velocitynumber_{\text{cons}} - 1, n}}, \quad &\populationletter = \velocitynumber_{\text{cons}}, \dots, \velocitynumber - 1. \label{eq:LeavesCollision}
\end{align}
This is significantly cheaper than \eqref{eq:reconstructedCollision} because there is no need to reconstruct a piece-wise constant representation of the solution on the full finest level.
Moreover, in the case where the equilibria are linear, we are still able to prove Proposition \ref{proposition:ErrorEstimate} and we might argue that in practice it is still verified for non-linear cases.
We shall validate this claim with simulations and provide an \emph{ad-hoc} pathological example where Proposition \ref{proposition:ErrorEstimate} does not hold, with full discussion in the Supplementary material.

For the stream phase, even if we reconstruct at the finest level, the computation can be done at minimal expenses because we are capable of passing from \eqref{eq:StreamUncompacted} to \eqref{eq:AdaptiveStreamPhase} by linearity.
Using cheaper reconstruction operators as hinted by Remark \ref{rem:TooCherStream} cannot yield the control by Proposition \ref{proposition:ErrorEstimate} and we have verified that it frequently generates low-quality results. This is the subject of a future contribution.

The algorithms are sequentially implemented in \texttt{C++} using a code called \texttt{SAMURAI}\footnote{Code, test cases and documentation available at \url{https://github.com/hpc-maths/samurai}.} (\textbf{S}tructured \textbf{A}daptive mesh and \textbf{MU}lti-\textbf{R}esolution based on \textbf{A}lgebra of \textbf{I}ntervals) which is currently under development and that can handle general problems involving dynamically refined meshes (both MR and AMR).
The central features of \texttt{SAMURAI} are its data structure based on intervals of contiguous cells along each axis and an ensemble of set operations to quickly and easily perform inter-level operations.

\section{Verifications}\label{sec:Verifications}

In this Section, we concentrate on two main aspects, namely:
\begin{itemize}
    \item The fulfillment of the theoretical estimate by Proposition \ref{proposition:ErrorEstimate}.\footnote{Even when we are not able to verify the continuity property of the reference scheme.}
        The errors are measured on the conserved moments.
        In particular, we look at:
        \begin{align*}
        E^{\populationletter, n} &\definitionequality \lVert \vectorial{\average{M}}^{\populationletter, \text{ex}}(t^n) - \vectorial{\average{M}}^{\populationletter, n} \rVert_{\ell^1}, \quad  e^{\populationletter, n} \definitionequality \lVert \vectorial{\average{M}}^{\populationletter, \text{ex}}(t^n) - \reconstructed{\vectorial{\average{\momentumadaptive}}}{\populationletter, n}{\maxlevel} \rVert_{\ell^1}, \\
        \delta^{\populationletter, n} &\definitionequality \lVert \vectorial{\average{M}}^{\populationletter, n} - \reconstructed{\vectorial{\average{\momentumadaptive}}}{\populationletter, n}{\maxlevel} \rVert_{\ell^1},
        \end{align*}
        for $\populationletter = 0, \dots, \velocitynumber_{\text{cons}} - 1$, which are respectively the error of the reference method against the exact solution \additionRefOne{(called ``reference discretization error'')}\additionRefOne{, the error of the adaptive method against the exact solution (called ``adaptive discretization error'')} and the difference between the adaptive solution and the reference solution \additionRefOne{(called ``perturbation error'')}.
        \additionRefOne{As seen in Proposition \ref{proposition:ErrorEstimate}, the perturbation error $\delta^{\populationletter, n} \to 0$ as $\thresholdletter \to 0$.}
        \additionRefOne{By the triangle inequality, we have that $e^{\populationletter, n} \leq E^{\populationletter, n} + \delta^{\populationletter, n}$. An important aspect once utilizing multiresolution, linked with the choice of $\thresholdletter$, is not to perturb the reference discretization error due to the perturbation error, \emph{e.g.} having $\delta^{\populationletter, n} \ll E^{\populationletter, n}$. This is independent of the fact that the reference scheme is convergent (and many lattice Boltzmann schemes are not), namely $E^{\populationletter, n} \to 0$ as $\maxlevel \to +\infty$.}
        \additionRefOne{Clearly, for convergent schemes, if the user increases $\maxlevel$, the threshold parameter $\epsilon$ has to be decreased accordingly in order to avoid interference with the convergence of the scheme, thus to have $\delta^{\populationletter, n} \ll E^{\populationletter, n}$ entailing $e^{\populationletter, n} \approx E^{\populationletter, n}$.}
        
        \additionRefTwo{We measure the differences on the conserved moments because they are ultimately the $\velocitynumber_{\text{cons}}$ quantities we are interested in and for which the exact solution is available.}
        \item The gain in terms of computational time induced by the use of multiresolution. 
        In this work, we use the compression factor \additionRefOne{(at final time)}, which is given by $\additionRefOne{\text{CR}^{n} \definitionequality } 100 \times (1-\ratioobl{\#{\physicalleaves{\Lambda^{n}}}}{N_{\maxlevel}})$, as a measure of computational efficiency, knowing that the real one is strongly dependent on the implementation and data structure and will be studied in future works.
        \additionRefOne{Observe that we also can use the time-average compression factor given by $ \text{ACR}^{N} \definitionequality 100 \times(1- 1/N \sum_{n=1}^N \#(S(\Lambda^n))/N_{\overline J} )$. In what we did, this metric is generally bounded from below by the compression factor at the final time for the following reason. We mostly start from solutions with shocks, where very high compression rates are achieved. Eventually, we obtain travelling shocks, contact discontinuities and rarefactions fans, thus having to put more and more cells and worsening the compression rate. Thus, the compression rate at the final time clearly bounds the average compression rate from below.}

\end{itemize}
Unless otherwise stated, the test are carried using the ``leaves collision''.
An exception to this rule is presented in details in the Supplementary material.
\additionFree{In the manuscript, we consider $\predictionstencildepth = 1$ (the number of neighbors considered by the prediction operator), $\thresholdletter = 10^{-4}$ (except when this parameter is varied or otherwise said), $\minlevel = 2$ and $\maxlevel = 9$. This value for $\thresholdletter$ guarantees, for any considered test with reference scheme at level $\maxlevel = 9$, to achieve $\delta^{h, n} \ll E^{h, n}$.}
\additionRefOne{Notice that $\maxlevel$ determines which reference scheme we are relying on and building our adaptation strategy. On the other hand, $\minlevel$ such that $0 \leq \minlevel \leq \maxlevel$ can be chosen freely by the user and determines the number of potential levels in the adaptive mesh.}
\additionRefTwo{Some test cases are repeated using a larger prediction stencil $\predictionstencildepth = 2$ in the Supplementary material.}
\additionRefOne{Larger prediction stencils entail larger costs of the multiresolution analysis and thus a larger overhead. However, since the decay of the details is linked to $\predictionstencildepth$ \emph{via} \eqref{eq:DetailsDecayEstimate} through $\mu = 2\predictionstencildepth + 1$, using $\predictionstencildepth = 2, 3, \dots$ for very smooth solutions can be beneficial to achieve very high compression factors. On the other hand, for solutions with shocks, the details do not decay with the level $\levelletter$ whatever $\predictionstencildepth$ is. It is therefore not advisable to use large $\predictionstencildepth$ in this situation, because this does not yield important gains in the mesh compression compared to the larger overhead.}

\additionRefOne{In the numerical simulations, many different schemes with several choices of relaxation parameters are considered. The aim is to showcase the generality of our approach. We do not focus on the choice of relaxation parameters to obtain the stability of the scheme with a good compromise between spurious oscillations and numerical diffusion since this is a huge subject on its own.}

\subsection{\scheme{1}{2} for a scalar conservation law: advection and Burgers equations}
\subsubsection{The problem and the scheme}

We aim at approximating the weak entropic solution (see Serre \cite{serre1999}) of the initial-value problem:
\begin{equation}\label{eq:ScalarConservation}
    \begin{cases}
        \partial_t \rho + \partial_x (\varphi(\rho)) = 0, \qquad t \in [0, T], \quad &x \in \reals, \\
        \rho(t = 0, x) = \rho_0(x), \qquad &x \in \reals.
    \end{cases}
\end{equation}
 with $\varphi \in C^{\infty} (\reals)$ a flux and $\rho_0 \in L^{\infty} (\reals)$.
This problem is the advection equation with constant velocity for $\varphi(\xi) = c \xi$ with velocity $c \in \reals$ and the inviscid Burgers equation for $\varphi(\xi) = \ratioobl{\xi^2}{2}$.
The \scheme{1}{2} scheme is obtained by selecting $\velocitynumber = 2$ and $\velocitynumber_{\text{cons}} = 1$ with velocities $\velocityletter^0 = \latticevelocity$, $\velocityletter^1 = -\latticevelocity$ and change of basis
\begin{equation*}
    \operatorial{M} = \adaptiveroundbrackets{\begin{matrix}
                                              1 & 1 \\
                                              \latticevelocity & - \latticevelocity
                                             \end{matrix}}.
\end{equation*}
With the theory of equivalent equations \cite{dubois2009}, Graille \cite{graille2014} has shown that the equivalent equation for this scheme is \eqref{eq:ScalarConservation} up to first order in $\Delta t$ upon selecting $M^{1, \text{eq}} = \varphi(\average{M}^0)$.

\begin{example}
    In the case of advection equation with $\latticevelocity \geq c > 0$, we have an explicit expression for the optimal continuity constant of the scheme for the $\ell^1$ norm, namely $C_L = 1$ if $s^1 \leq \ratioobl{2}{(1 + \ratioobl{c}{\latticevelocity})}$  or $C_L = s^1 \adaptiveroundbrackets{1 + \ratioobl{c}{\latticevelocity}} - 1$ otherwise.
\end{example}

\subsubsection{Results}

        \begin{table}[h]\caption{\label{tab:TestCases}Test cases for one scalar conservation law with choice of flux, initial datum, expected regularity of the solution, choice of the regularity parameter $\overline{\mu}$ and final time of the simulation.}
        \begin{footnotesize}
        \begin{tabular}{|c|c|c|c|c|c|}
        \hline
        Flux $\varphi$ & Initial datum $\rho_0$ & Type of solution & $\nu$ & $T$ & Test \\
         \hline 
         \hline
         \multirow{2}{*}{$\varphi(u) = \dfrac{3}{4}u$} & $\rho_0(x) = e^{-20x^2}$ & Strong $C^{\infty}$ & $\infty$ & 0.4 & \textrm{I} \\
         & $\rho_0(x) = \chi_{|x| \leq 1/2}(x)$ & Weak $L^{\infty}$ & 0 & 0.4 & \textrm{II} \\
         \hline
         \multirow{3}{*}{$\varphi(u) = \dfrac{u^2}{2}$} & $\rho_0(x) = (1+\text{tanh}(100x))/2$ & Strong $C^{\infty}$ & $\infty$ & 0.4 & \textrm{III} \\
          & $\rho_0(x) = \chi_{|x| \leq 1/2}(x)$ & Weak $L^{\infty}$ & 0 & 0.7 & \textrm{IV} \\
          & $\rho_0(x) = (1+x)\chi_{x < 0}(x) + (1-x)\chi_{x \geq 0}(x)$ & Weak $L^{\infty}$ & 0 & 1.3 & \textrm{V} \\
          \hline
        \end{tabular}
        \end{footnotesize}
    \end{table}

                \begin{figure}
                \begin{center}
                    (I) \\
                    \includegraphics[width=1.0\textwidth]{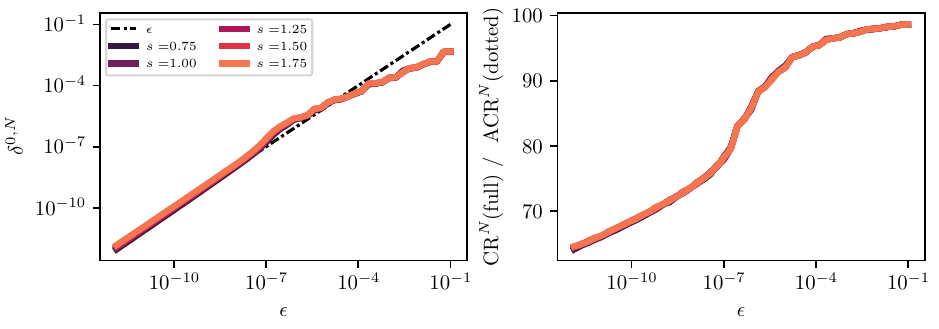}
                    \vspace{-0.05cm}(II) \\
                    \includegraphics[width=1.0\textwidth]{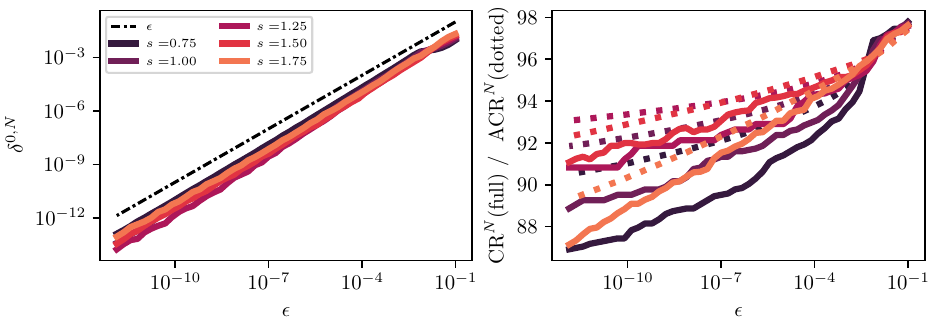}
                    \vspace{-0.05cm}(III) \\
                    \includegraphics[width=1.0\textwidth]{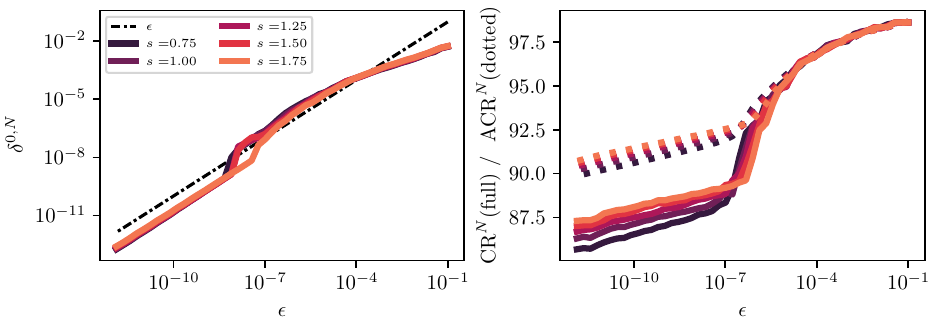}
                    \vspace{-0.05cm}(IV) \\
                    \includegraphics[width=1.0\textwidth]{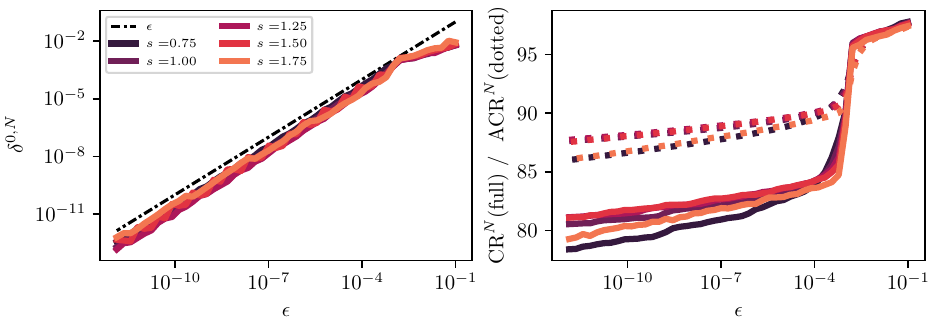}
                \end{center}\caption{\label{fig:epsilondifferentests}Behavior of $\delta^{0, N}$ as function of $\thresholdletter$ (left) and compression factors (full line for $\text{CR}^N$ and dotted line for $\text{ACR}^N$) at the final time as function of $\thresholdletter$ (right), for test (from top to bottom) \textrm{I}, \textrm{II}, \textrm{III} and \textrm{IV}. The dot-dashed line gives the reference $\delta^{0, N} = \thresholdletter$. We call $s = s^1$.}
            \end{figure}

            \begin{figure}
                \begin{center}
                    \includegraphics[width=1.0\textwidth]{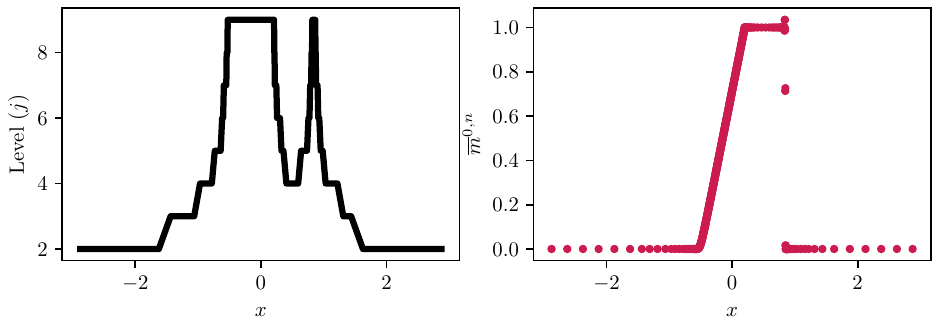}
                \end{center}\caption{\label{fig:solutionDoubleRiemann}Example of solution of the \scheme{1}{2} for the Test \textrm{IV}, considering $n = 358$, $s^1 = 1.5$ and $\thresholdletter = 10^{-4}$. On the left, levels of the computational mesh. On the right, solution on the leaves of the tree.}
            \end{figure}


            \begin{table}\caption{\label{tab:RatioErrors}Value of the ratio $E^{0, N}/\delta^{0, N}$ at final for each test case for one scalar conservation law. The time variation of this quantities can be found in the Supplementary material.}
        \begin{footnotesize}
            \begin{center}
            \begin{tabular}{|c||c|c|c|c|c|}
            \hline
             & \multicolumn{5}{c|}{$E^{0, N}/\delta^{0, N}$} \\
             \hline
            $s$ & \textrm{I} & \textrm{II} & \textrm{III} & \textrm{IV} & \textrm{V}\\
            \hline
0.75 &	9.97e+01 &	1.86e+03 &	5.93e+01 &	3.50e+02 &	8.78e+02 \\
1.00 &	5.94e+01 &	2.31e+03 &	3.71e+01 &	3.41e+02 &	1.01e+03 \\
1.25 &	3.52e+01 &	2.62e+03 &	2.29e+01 &	3.93e+02 &	9.89e+02 \\
1.50 &	1.94e+01 &	2.44e+03 &	1.31e+01 &	9.72e+01 &	1.05e+03 \\
1.75 &	8.34e+00 &	1.21e+03 &	5.71e+00 &	2.90e+02 &	1.14e+03 \\
            \hline
            \end{tabular}
            \end{center}
        \end{footnotesize}
    \end{table}

    For this test case, we consider $\Omega = [-3, 3]$.
    Concerning the lattice Boltzmann scheme, we fix $\latticevelocity = 1$ and we vary the relaxation parameter.
    The tests we perform are resumed on Table \ref{tab:TestCases}.
    \additionRefOne{Observe that the choice for the regularity guess $\overline{\mu}$ has been done and shall be done according to the expected smoothness of the solution, namely, it should be equal to the number of expected bounded derivatives.}
    In order not to overcharge the paper with plots, we just provide the values for the ratio $E^{0, N}/\delta^{0, N}$ at final time $T$ on Table \ref{tab:RatioErrors}. The time evolution of $\delta^{0, n}$ as well as that of \additionRefOne{$E^{0, n}/\delta^{0, n}$}, $E^{0, n}/e^{0, n}$ \additionRefOne{and $\text{CR}^n$} can be found in the Supplementary material.
    The evolution of the perturbation error of the adaptive MR-LBM method and the compression factor as function of $\thresholdletter$ are given on Figure \ref{fig:epsilondifferentests}, except for test number V, which is discussed in more detail in the Supplementary material.    
    We formulate the following remarks:
    \begin{enumerate}
            \item[\textrm{I}.] 
            We observe that with this choice of $\thresholdletter$ we successfully keep the perturbation error by the adaptive MR-LBM scheme $\delta^{0, n}$ about 10-100 times smaller than the reference discretization error $E^{0, n}$ at the chosen level $\maxlevel$, with important compression rates around $95$\% for the chosen $\thresholdletter$.
            We remark the fairly correct linear behavior in terms of $\thresholdletter$.
            \additionRefOne{The compression factor $\text{CR}^N$ and the average compression factor $\text{ACR}^N$ coincide because the solution retains the same smoothness in time and is simply transported (plus the numerical diffusion).}
            We have verified that the perturbation error increases linearly\footnote{\label{note:supplmaterial}See Supplementary material.} in time even when we can only prove an exponential bound\footnote{\emph{I.e.} $s > \ratioobl{8}{7}$.} by Proposition \ref{proposition:ErrorEstimate}.
            
            \item[\textrm{II}.] The \additionRefOne{perturbation} error of the adaptive MR-LBM method is about three orders of magnitude smaller than the reference discretization error, again for the selected $\maxlevel$.
            Due to the presence of large \emph{plateaux}, the compression factor is really interesting for a large range of $\thresholdletter$, being always over $90$\%.
            \additionRefOne{We see that $\text{ACR}^N$ is larger than $\text{CR}^N$ arguably because of the numerical diffusion which accumulates in time and smears the shock.}
            The trend of $\delta^{0, N}$ as function of $\thresholdletter$ agrees with the theory and can be bound linearly in time.
            
            \item[\textrm{III}.] Again, we observe that this choice of $\thresholdletter$ guarantees perturbation errors which are between 5 and 50 times smaller than the \additionRefOne{discretization} error of the reference method, still preserving excellent compression rates.
            \additionRefOne{We again have $\text{ACR}^N > \text{CR}^N$ because of the formation of a rarefaction fan as the simulation goes on.}
            The behavior as $\thresholdletter$ tends to zero is respected and the expected linear temporal trend is obtained.
            
            \additionRefOne{The attentive reader could have observed the following fact: the compression rates tend to stagnate as $\thresholdletter \to 0$. The reason for that is the following (and applies to any context in which the situation shall happen, also in the sequel): consider the typical solution of most of the problems we consider, where only shocks (and contact discontinuities) and rarefaction fans are present. Elsewhere, the solution is essentially flat. Start from a very large threshold $\thresholdletter$: multiresolution does not put cells at the finest level of resolution $\maxlevel$ because the threshold is really large. Then decrease $\thresholdletter$ little by little: the finest resolution is reached on the shock and in the less smooth zones of the rarefaction fans. By continuing decreasing $\thresholdletter$, the fans are also refined (especially if here the solution is highly non-linear). Nevertheless, at some time, the finest level $\maxlevel$ is reached everywhere where the solution is non flat (shocks and fans) and eventually (for smaller $\thresholdletter$) there is not so much room for improving the quality of the reconstruction by refining elsewhere, because here the solution is totally flat (and indeed the details are perfectly equal to zero).
  This is why the compression rate (almost) stagnates. Multiresolution can still diminish the error as expected by adding very few cells thus with very little modifications of the compression rates. Of course, one expects $\text{CR}^N, \text{ACR}^N \to 0$ as $\epsilon \to 0$, but in this case $\epsilon$ should become really small, allegedly below the machine epsilon to observe the convergence after the stagnation.}
        
            \item[\textrm{IV}.] For illustrative purposes the weak solution of the problem is shown on Figure \ref{fig:solutionDoubleRiemann}). 
            The adaptive method largely beats the traditional method by three orders of magnitude, with less efficient compression compared to (\textrm{II}) due to the formation of a rarefaction fan\footnote{This rarefaction is straight-shaped but multiresolution refines at the extremal kinks of the slope. Moreover the \scheme{1}{2} creates a stair-shaped rarefaction, which triggers refinement.}, \additionRefOne{which is again the cause of  $\text{ACR}^N > \text{CR}^N$.}
            The estimate in $\thresholdletter$ is sharply met and the perturbation error increases linearly in time for every choice of $s^1$.\footnote{Sometimes with strong oscillations due to the oscillations of the scheme.}
            
            \item[\textrm{V}.] The outcome of this test is presented and fully discussed in the Supplementary material and provides a pathological example where the reconstructed collision is needed to correctly retrieve the theoretical estimates \additionRefOne{on the perturbation error}.
        \end{enumerate}
        
        Overall, we can conclude that the adaptive MR-LBM for a scalar conservation law guarantees an error control by a threshold $\epsilon$ and succeeds in keeping the perturbation error $\delta^{0, n}$ way smaller than the discretization error $E^{0, n}$ \additionRefOne{of the reference scheme }(see Table \ref{tab:RatioErrors}) especially when weak solutions are involved, \additionRefOne{for the selected maximum level $\maxlevel$}.
        The ``leaves collision'' does not impact these characteristics except in a specifically designed pathological case.

    \subsection{\scheme{1}{3} and \scheme{1}{5} for two conservation laws: the shallow water system}
    
        \subsubsection{The problem and the scheme}
            We aim at approximating the weak entropic solution of the shallow water system, where $h$ represent the height of a fluid and $u$ is its horizontal velocity:
            \begin{equation}\label{eq:ShallowWaters}
                \begin{cases}
                    \partial_t h + \partial_x (hu) = 0, \qquad t \in [0, T], \quad &x \in \reals, \\
                    \partial_t (hu) + \partial_x (hu^2 + gh^2/2) = 0,  \qquad t \in [0, T], \quad &x \in \reals, \\
                    h(t=0, x) = h_0(x), \qquad &x \in \reals, \\
                    u(t=0, x) = u_0(x), \qquad &x \in \reals, 
                \end{cases}
            \end{equation}
 where $g > 0$ is the gravitational acceleration exerted on the fluid and $h_0, u_0 \in L^{\infty}(\reals)$.
    Two possible lattice Boltzmann schemes with two conserved moments are:
   \begin{itemize}
            \item \scheme{1}{3}, obtained selecting $\velocitynumber= 3$ and $\velocitynumber_{\text{cons}} = 2$ with discrete velocities $\velocityletter^0 = 0$, $\velocityletter^1 = \latticevelocity$, $\velocityletter^2 = -\latticevelocity$ and the change of basis:
            \begin{equation*}
               \operatorial{M} = \adaptiveroundbrackets{\begin{matrix}
                                                                                    1 & 1 & 1 \\
                                                                                    0 & \latticevelocity & -\latticevelocity \\
                                                                                    0 & \latticevelocity^2 & \latticevelocity^2
                                                                                \end{matrix}}.
            \end{equation*}
            Selecting $M^{2, \text{eq}} = (\average{M}^1)^2/\average{M}^0 + g(\average{M}^0)^2/2$ the scheme is consistent up to first order in $\Delta t$ with \eqref{eq:ShallowWaters}.
            \item \scheme{1}{5}, obtained taking $\velocitynumber = 5$ and $\velocitynumber_{\text{cons}} = 2$ with the choice of velocities $\velocityletter^0 = 0$, $\velocityletter^1 = \latticevelocity$, $ \velocityletter^2 = -\latticevelocity$, $\velocityletter^3 = 2\latticevelocity$, $ \velocityletter^4 = -2\latticevelocity$, along with the matrix:
            \begin{equation*}
                \operatorial{M} = \adaptiveroundbrackets{\begin{matrix}
                                                                             1 & 1 & 1 & 1 & 1 \\
                                                                             0 & \latticevelocity & -\latticevelocity & 2 \latticevelocity & -2 \latticevelocity \\
                                                                             0 & \latticevelocity^2 & \latticevelocity^2 & 4\latticevelocity^2 & 4\latticevelocity^2 \\
                                                                             0 & \latticevelocity^3 & -\latticevelocity^3 & 8\latticevelocity^3 & -8\latticevelocity^3 \\
                                                                             0 & \latticevelocity^4 & \latticevelocity^4 & 16\latticevelocity^4 & 16\latticevelocity^4 \\
                                                                            \end{matrix}}.
            \end{equation*}
            We select the equilibri in the following way:
            \begin{equation*}
                M^{2, \text{eq}} = \frac{(\average{M}^1)^2}{\average{M}^0} + \frac{g}{2} (\average{M}^0)^2, \qquad M^{3, \text{eq}} = \alpha \latticevelocity^2 \average{M}^1, \qquad M^{4, \text{eq}} = \beta \latticevelocity^2 \average{M}^{2, \text{eq}},
            \end{equation*}
            where $\alpha$ and $\beta$ are real parameters to be set in order to keep the scheme stable. With this choice the equivalent equations are consistent with \eqref{eq:ShallowWaters} up to first order being close to those of the \scheme{1}{3} scheme.
        \end{itemize}
        
\subsubsection{Results}
              \begin{figure}
                \begin{center}
                    \includegraphics[width=1.0\textwidth]{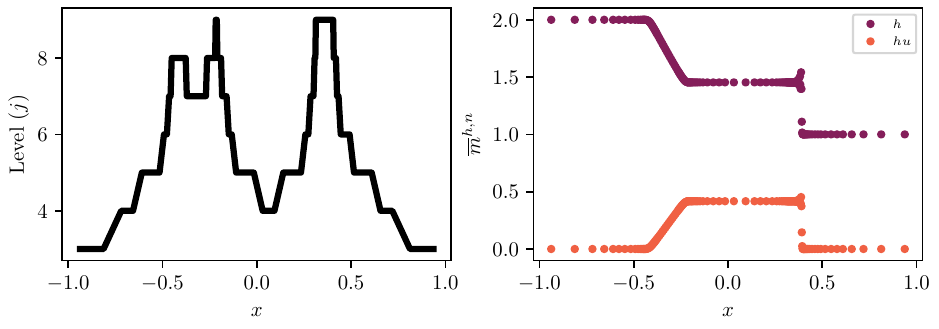}
                \end{center}\caption{\label{fig:plotSolutionSW}Example of solution of the \scheme{1}{5} for the shallow water problem with $n = 300$, $s^2 = 1.6$ and $\thresholdletter = 10^{-4}$. On the left, levels of the computational mesh. On the right, solution on the leaves of the tree.}
            \end{figure}

            \begin{figure}
                \begin{center}
                    \includegraphics[width=1.0\textwidth]{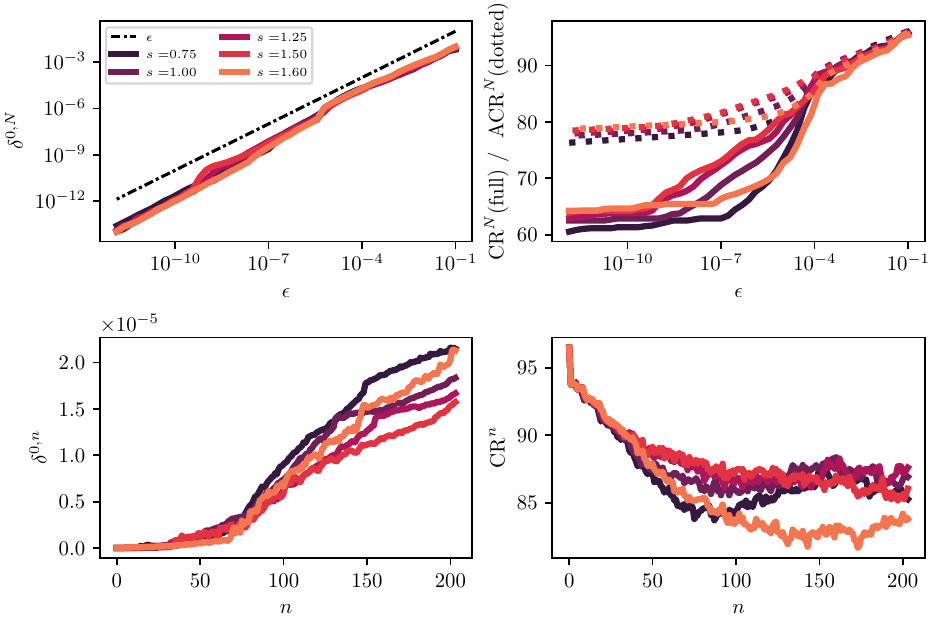}
                \end{center}\caption{\label{fig:d1q3_s_d}\scheme{1}{3} for the shallow water system with Riemann initial datum. The dot-dashed line gives the reference $\delta^{0, N} = \thresholdletter$. For the sake of avoiding redundancy, we present only one moment. \additionRefOne{$s$ is a shorthand for $s^2$.}}
            \end{figure}
            
             \begin{figure}
                \begin{center}
                    \includegraphics[width=1.0\textwidth]{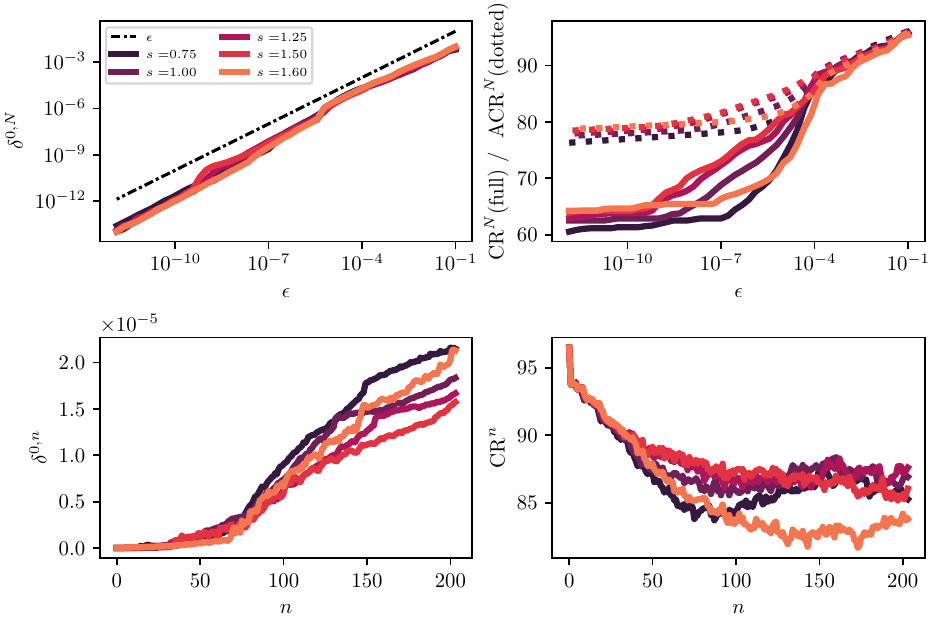}
                \end{center}\caption{\label{fig:d1q5_s_d}\scheme{1}{5} for the shallow water system with Riemann initial datum. The dot-dashed line gives the reference $\delta^{0, N} = \thresholdletter$. For the sake of avoiding redundancy, we present only one moment. \additionRefOne{$s$ is a shorthand for $s^2$.}}
            \end{figure} 
            
            As initial datum, we consider the Riemann problem given by $(h_0, u_0) (x) =  (2, 0)\chi_{x < 0}(x) + (1, 0)\chi_{x>0}(x)$, with a lattice velocity $\latticevelocity = 2$, a final time $T = 0.2$ and a domain $\Omega=[-1, 1]$.
            \additionFree{The gravitational acceleration is $g = 1$.}
            The result is shown in Figure \ref{fig:d1q3_s_d}: for both the conserved moments, the behavior of the perturbation error in time is supra-linear, being very small at the very beginning because the method adds enough security cells around the shock and information propagates relatively slowly.
            Moreover, we remark that the \additionRefOne{perturbation} error is larger for smaller $s^2$ due to the larger diffusivity of the numerical scheme.
            The perturbation error is between four and six orders of magnitude smaller than the \additionRefOne{discretization} error of the reference method, reaching very interesting compression factors.
            The estimates \additionRefOne{for $\delta^{h, N}$ for $h = 0, 1$} in terms of $\thresholdletter$ are correctly followed.
            \additionRefOne{We observe the typical inequality $\text{ACR}^N > \text{CR}^N$.}
            
            For the \scheme{1}{5}, we use exactly the same setting except for taking $\alpha = \beta = 1$ and setting $s^3 = s^4 = 1$.
            We obtain what is shown in Figures \ref{fig:plotSolutionSW} and \ref{fig:d1q5_s_d}.\footnote{We are limited to $s^2 = 1.6$ due to stability issues which are inherent to the reference scheme.}
            The time behavior of the perturbation error is again supra-linear and now the difference between different relaxation parameters is less evident.
            The ratio with the \additionRefOne{discretization} error of the reference scheme is between $10^4$ and $10^6$.
            The bound \additionRefOne{of $\delta^{h, N}$ for $h = 0, 1$} in $\thresholdletter$ is very well fulfilled.
            This example shows that our adaptive strategy works really well even for schemes with an extended advection stencil, namely with $\sigma = 2$.

\subsection{\schemevec{1}{2}{3} for the Euler system}

We consider the full Euler system
\begin{equation}\label{eq:EulerSystem}
    \begin{cases}
        \partial_t \rho + \partial_x(\rho u) = 0, \qquad t \in [0, T], \quad &x \in \reals, \\
        \partial_t (\rho u) + \partial_x (\rho u^2 + p) = 0,  \qquad t \in [0, T], \quad &x \in \reals, \\
        \partial_t E + \partial_x (Eu + pu) = 0,  \qquad t \in [0, T], \quad &x \in \reals, \\
        \rho(t=0, x) = \rho_0(x), \qquad &x\in \reals,\\
        u(t=0, x) = u_0(x), \qquad &x\in \reals,\\
        E(t=0, x) = E_0(x), \qquad &x\in \reals,
    \end{cases}
\end{equation}
        where $\rho$ is the density, $u$ the velocity of the flow, $p$ the pressure and $E$ the total energy.
        The pressure and the energy are linked by the pressure law $E = \ratioobl{\rho u^2}{2} + \ratioobl{p}{(\gamma_{\text{gas}} - 1)}$.
        For this work, we consider the Sod shock problem, choosing $\gamma_{\text{gas}} = 1.4$ and considering the Riemann initial datum given by:
        \begin{equation*}
            (\rho_0, u_0, E_0)(x) = (1.000, 0.000, 2.500) \chi_{x<0}(x) + (0.125, 0.000, 0.250) \chi_{x > 0}(x),
        \end{equation*}
        generating a solution with a left-moving rarefaction, a right-moving contact discontinuity and a right-moving shock.
        We employ a vectorial scheme \cite{graille2014, dubois2014} rather than a scalar one for it adds the necessary numerical diffusion, enhancing stability and it makes easy to conserve $E$ without further manipulation.
        The scheme is the juxtaposition of three \scheme{1}{2} for the quantities $\rho$, $\rho u$ and $E$, coupled through their equilibri:
        \begin{align*}
            M^{1, \text{eq}} = \average{M}^2, \qquad M^{3, \text{eq}} &= \adaptiveroundbrackets{\ratio{3}{2} - \ratio{\gamma_{\text{gas}} }{2}} \ratio{(\average{M}^2)^2}{\average{M}^0} + (\gamma_{\text{gas}}  - 1) \average{M}^4, \\
            M^{5, \text{eq}} &= \gamma_{\text{gas}}  \ratio{ \average{M} ^4  \average{M}^2}{\average{M}^0} + \ratio{1-\gamma_{\text{gas}} }{2}\ratio{(\average{M}^2)^3}{(\average{M}^0)^2}.
        \end{align*}
        This scheme is consistent up to first order with \eqref{eq:EulerSystem} as shown by Graille \cite{graille2014}.

        \subsubsection{Results}
  
            \begin{figure}
                \begin{center}
                    \includegraphics[width=1.0\textwidth]{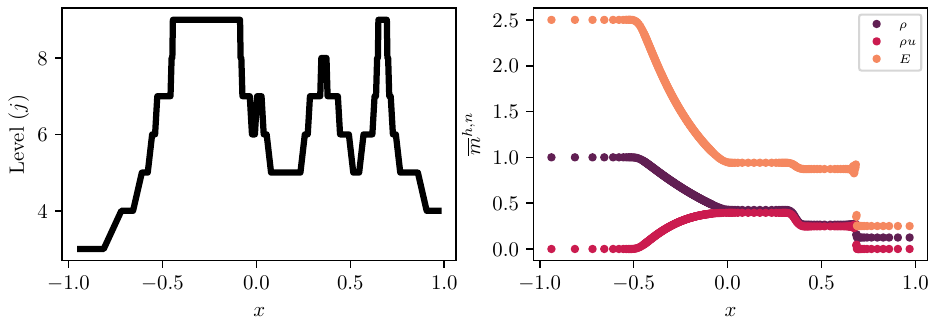}
                \end{center}\caption{\label{fig:plotSolutionSod}Example of solution of the vectorial \scheme{1}{2} for the Sod problem with $n = 600$, $s = 1.75$ and $\thresholdletter = 10^{-3}$. On the left, levels of the computational mesh. On the right, solution on the leaves of the tree.}
            \end{figure}
                    
            \begin{figure}
                \begin{center}
                    \includegraphics[width=1.0\textwidth]{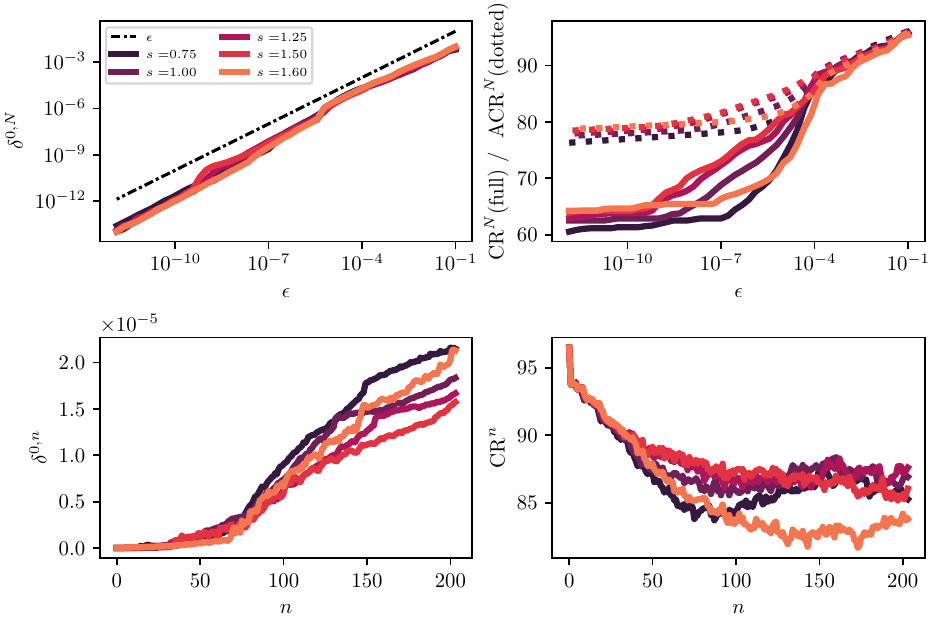}
                \end{center}\caption{\label{fig:d1q2_3_s_d}Vectorial \scheme{1}{2} for the Sod problem. The dot-dashed line gives the reference $\delta^{2, N} = \thresholdletter$. For the sake of compactness, we show only one moment.}
            \end{figure}
        
        We consider a domain $\Omega = [-1, 1]$ and all the other parameters as in the previous examples, except the lattice velocity taken to be $\latticevelocity = 3$ and the final time $T = 0.4$. 
        All the relaxation parameters are taken equal\additionRefOne{: call them $s$}.
        The result is given in Figures \ref{fig:plotSolutionSod} and \ref{fig:d1q2_3_s_d}: the perturbation error behaves fairly linearly in time for every choice of relaxation parameter and becomes smaller as $s$ approaches two, due to the reduced numerical diffusion.
        We are capable of keeping the perturbation error between three and four order of magnitudes smaller than to the \additionRefOne{discretization} error of the reference scheme for each of the conserved moments, \additionRefOne{for the chosen resolution $\maxlevel$}.
        The behavior in $\thresholdletter$ is respected.        
        This shows that our strategy is well suited to handle the simulation of systems of conservation laws using vectorial schemes.

\section{Conclusions}\label{sec:Conclusions}

    
    In this paper, we have presented a class of new fully adaptive lattice Boltzmann schemes based on multiresolution to perform the adaptation of the spatial grid with error control.
    To the best of our knowledge, no previous research has been conducted to couple multiresolution and LBM methods.
    The most important features are that there is no need to devise \emph{ad-hoc} refinement/coarsening criteria: mesh adaptation is naturally handled using multiresolution by analyzing the regularity of the solution.
    Therefore, no previous knowledge of the solution\footnote{Other than the regularity guess $\overline{\mu}$, which can be set to a small value for precaution.} is needed because the numerical mesh is automatically evolved.
    Eventually, under reasonable assumptions on the reference scheme, we are able to prove precise error controls on the perturbation error introduced by the non-uniform mesh, which are driven by a single adjustable tolerance $\thresholdletter$.
    The numerical method has been extensively tested, showing that the theoretical predictions are fully met, even for settings for which we expect less predictable behaviors.
    We have shown that, by tuning $\thresholdletter$ according to the desired precision, one is capable of keeping the perturbation error several orders of magnitude below the discretization error of the reference method with respect to the exact solution, still achieving excellent compression factors.
    We have also demonstrated that the optimized ``leaves collision'' is an efficient alternative to the ``reconstructed collision'', except in pathological cases.
    The question on how the choice of prediction operator could modify the physics approximated by the MR-LBM adaptive scheme will be the object of a forthcoming contribution.

    
    The major improvement our method needs to undergo is its generalization to the multi-dimensional framework (spatial dimension $d$). We provide answers to this question in a companion contribution \cite{bellotti2021}.
    The following points have been taken into consideration:
    \begin{itemize}
        \item The projection operator is straightforwardly generalized as a mean on the children. The prediction operator is constructed by tensor product as hinted by Bihari and Harten \cite{bihari1997}.
        \item The decay estimates for the details \eqref{eq:DetailsDecayEstimate} are still valid without having to adjust them with $d$.
        Consequently \eqref{eq:DetailRatio} remains valid.
        However, one shall cope with the fact that details of two siblings no longer have the same modulus.
        \item One must modify the choice of $\thresholdletter_{\levelletter}$ according to $d$, as the number of elements in a tree is now bounded by $2^{d\maxlevel}$. 
        We consider $\thresholdletter_{\levelletter}= 2^{d (\levelletter - \maxlevel)}\thresholdletter$.
        Hence, we need to slightly modify \eqref{eq:HartenLossOfRegularity} which becomes $|\average{\detailsadaptive}_{\levelletter, \indexletter}^{\populationletter, n}| \geq 2^{\overline{\mu} + d} \thresholdletter_{\levelletter}$.
        \item The stream phase given by \eqref{eq:AdaptiveStreamPhase} and the way of recovering it remain essentially the same.
    \end{itemize}
    In \cite{bellotti2021}, we employ the MR-LBM adaptive scheme to simulate both hyperbolic (Euler) and parabolic (incompressible Navier-Stokes) systems, because the accuracy of our reconstruction is enough to correctly cope with the physics of such systems.
    
    Finally, the optimisation of the implementation is a crucial subject when dealing with multidimensional problems.
    In this work, we have restricted purposefully the measure of the computational gain with respect to the uniform mesh by merely looking at the compression factor\additionRefOne{s}.
    This is far from realistic if the implementation does not perform the operations involved in multiresolution in a clever way or if the problem is too small to observe a real gain.
    We believe that the choice of the underlying data structure has a huge impact on this matter: we are currently developing the library \texttt{SAMURAI} with the purpose of providing an innovative interval-based data structure to enhance performances and simplify the parallelization  of the whole process. This is the subject of our current research.
    
\section{Acknowledgements}
    The authors deeply thanks Laurent Séries for fruitful discussions on multiresolution. \additionFree{They also thank the two anonymous referees for the useful remarks and suggestions.}
    Thomas Bellotti is supported by a PhD funding (year 2019) from the Ecole polytechnique.

\bibliographystyle{siam}
\bibliography{bibliography.bib}

\end{document}